\documentclass[leqno,12pt]{article}
\usepackage{amssymb, amsfonts, amsmath}
\usepackage{xcolor}
\usepackage{hyperref}
\usepackage[numbers]{natbib}
\newcommand{\ra}{\rightarrow}
\newcommand{\Ra}{\Rightarrow}

\newcommand{\lb}{\langle}
\newcommand{\rb}{\rangle}

\newcommand{\ov}{\overline}
\newcommand{\sref}[1]{(\ref{#1})}

\newcommand{\bmat}{\left(\begin{array}}
\newcommand{\emat}{\end{array}\right)}
\newcommand{\beq}{\begin{equation}}
\newcommand{\eeq}{\end{equation}}

\newcommand{\pf}{\noindent{\em Proof.}}
\newcommand{\pff}{\noindent{\em Proof of}}

\newcommand{\Mod}[1]{\,(\mbox{\rm mod}~#1)}
\newcommand{\C}{\mathbb{C}}

\newcommand{\Q}{\mathbb{Q}}
\newcommand{\Z}{\mathbb{Z}}
\newcommand{\cA}{{\cal A}}

\newcommand{\cN}{{\cal N}}

\newcommand{\cP}{{\cal P}}

\newcommand{\al}{\alpha}
\newcommand{\be}{\beta}
\newcommand{\de}{\delta}
\newcommand{\De}{\Delta}
\newcommand{\ga}{\gamma}
\newcommand{\Ga}{\Gamma}
\newcommand{\ve}{\varepsilon}
\newcommand{\lm}{\lambda}
\newcommand{\vp}{\varphi}
\newcommand{\om}{\omega}
\newcommand{\Om}{\Omega}
\newcommand{\tth}{\theta}
\newcommand{\Th}{\Theta}
\newcommand{\fa}{\mathfrak{a}}
\newcommand{\fA}{\mathfrak{A}}
\newcommand{\fc}{\mathfrak{c}}
\newcommand{\fC}{\mathfrak{C}}
\newcommand{\fH}{\mathfrak{H}}
\newcommand{\fm}{\mathfrak{m}}
\newcommand{\fO}{\mathfrak{O}}
\newcommand{\fW}{\mathfrak{W}}
\newcommand{\Aut}{\operatorname{Aut}}
\newcommand{\Br}{\operatorname{Br}}
\newcommand{\Char}{\operatorname{char}}
\newcommand{\cl}{\operatorname{cl}}
\newcommand{\cont}{\operatorname{cont}}
\newcommand{\disc}{\operatorname{disc}}
\newcommand{\Div}{\operatorname{Div}}
\newcommand{\End}{\operatorname{End}}
\newcommand{\gen}{\operatorname{gen}}
\newcommand{\GL}{\operatorname{GL}}
\newcommand{\Hom}{\operatorname{Hom}}
\newcommand{\Ker}{\operatorname{Ker}}
\newcommand{\lcm}{\operatorname{lcm}}
\newcommand{\NS}{\operatorname{NS}}
\newcommand{\Pic}{\operatorname{Pic}}
\newcommand{\rank}{\operatorname{rank}}
\newcommand{\SL}{\operatorname{SL}}
\newcommand{\wt}{\operatorname{wt}}
\newcommand{\irr}{\operatorname{irr}}
\newcommand{\red}{\operatorname{red}}
\newcommand{\CM}{\operatorname{CM}}
\newcommand{\fun}{\operatorname{fun}}
\newcommand{\cpe}{\cP(A)^{\operatorname{ev}}}
\newcommand{\cpi}{\cP(A)^{\operatorname{irr}}}
\newcommand{\cpo}{\cP(A)^{\operatorname{odd}}}
\newcommand{\cpoi}{\cP(A)^{\operatorname{odd}, \operatorname{irr}}}
\newcommand{\cpr}{\cP(A)^{\operatorname{red}}}
\newcommand{\ocpi}{\ov\cP(A)^{\operatorname{irr}}}
\newcommand{\ocpoi}{\ov\cP(A)^{\operatorname{odd},\operatorname{irr}}}
\newcommand{\ocpr}{\ov\cP(A)^{\operatorname{red}}}
\newcommand{\thev}{\Theta_A^{\operatorname{ev}}}
\newcommand{\tho}{\Theta_A^{\operatorname{odd}}}
\newcommand{\leg}[2]{\left(\frac{#1}{#2}\right)}
\newtheorem {Theorem}{Theorem} %[section]
\newtheorem {thm}[Theorem]{Theorem} 

\newtheorem {prop}[Theorem]{Proposition}
\newtheorem {cor}[Theorem]{Corollary}
\newtheorem {lemma}[Theorem]{Lemma}

\newtheorem {rmk}[Theorem]{Remark}

\newtheorem {exam}[Theorem]{Example} 
\newcommand{\defn}{\noindent{\bf Definition.}}
\newcommand{\notn}{\noindent{\bf Notation.}}
\newcommand{\sps}{\vspace{3pt}}   %small space
\newcommand{\spm}{\vspace{6pt}} 
\newcommand{\noi}{\noindent}
\newcommand{\exprow}[1]{\renewcommand{\arraystretch}{#1}}
\newcommand{\expcol}[1]{\addtolength{\arraycolsep}{#1}} 
\newcommand{\harun}[1]{{\color{red}[Harun: #1]}}
\newcommand{\kani}[1]{{\color{blue}[Kani: #1]}}

\begin{document}
\vspace*{-30pt}

\begin{center}
\LARGE \textbf{The Number of Curves of Genus 2 with a Given Refined Humbert Invariant} \par
\vspace{20pt}

\large \textbf{Ernst Kani$^{\text{a}}$, Harun K\i r$^{\text{a}, 1, \ast}$} \par
\vspace{15pt}

\normalsize
$^{\text{a}}$ \textit{Department of Mathematics and Statistics, Queen's University, Kingston, Ontario, K7L 3N6, Canada. \\ E-mail addresses: \texttt{kani@mast.queensu.ca} (E. Kani), \texttt{harun.kir@ens-lyon.fr} (H. K\i r)}
\end{center}

\let\thefootnote\relax\footnotetext{$^1$ \textit{Current address:} ENS de Lyon, UMPA, UMR 5669, Lyon, France.}
\let\thefootnote\relax\footnotetext{$^\ast$ \textbf{Corresponding author.}}

\vspace{20pt}

 \begin{abstract}
Let $C/K$ be a curve of genus 2 over an algebraically closed field $K$.
Every such curve comes equipped with a canonical quadratic form $q_C$
called its refined Humbert invariant. In the case that the Jacobian $J_C$ of
$C$ is isogenous to the self-product $E \times E$ for an elliptic curve
$E/K$ with complex multiplication (CM), we provide an explicit formula
for the finite number $N_C$ of isomorphism classes of genus $2$ curves
$C'/K$ whose refined Humbert invariant is equivalent to $q_C$.  This formula implies that $N_C$ is unbounded for such curves, and that there are only finitely many isomorphism classes of such curves $C/K$ with a given value of $N_C$. 

A key step in our approach is a characterization of when $J_C$ is isogenous to a self product of a CM elliptic curve, formulated purely in terms of properties of the refined Humbert invariant $q_C$. We establish that an analogous characterization also holds for superspecial curves of genus 2. The paper concludes with explicit examples illustrating cases where a genus 2 curve is uniquely determined by the invariant $q_C$.
\end{abstract}

\iffalse
\begin{abstract}
    Let $C$ be a curve of genus 2 over an algebraically closed field. Every such curve comes equipped with a canonical quadratic form $q_C$ called its refined Humbert invariant. When the Jacobian of $C$ is isogenous to the self-product of a complex multiplication (CM) elliptic curve, we provide an explicit formula for the finite number $N_C$ of isomorphism classes of genus 2 curves whose refined Humbert invariant is equivalent to $q_C$. This formula implies that $N_C$ is unbounded for such curves, and that only finitely many isomorphism classes of such curves yield any given value of $N_C$.

A key step is a characterization of when the Jacobian is isogenous to a CM elliptic curve's self-product, formulated via properties of $q_C$. We establish an analogous characterization for superspecial genus 2 curves. The paper concludes with  examples where $C$ is uniquely determined by the invariant $q_C$.
\end{abstract}
\fi

\section{Introduction}
\label{s:intro}

If $C/K$ is a curve of genus $2$ over an algebraically closed field $K$, then it comes equipped with 
a canonical quadratic form $q_C$ called its \emph{refined Humbert invariant}; 
cf.\ \cite{K-EC}, \cite{K-JT}. This invariant is useful because many geometric properties of $C$ are   
reflected in arithmetic properties of the quadratic form $q_C$; 
cf.\ \cite{K-EC}, \cite{K-MJ}, \cite{K-SC}, \cite{K-CAS}, \cite{Ki_Aut} and \cite{GRV}. 

For example, the property that the Jacobian $J_C$ of $C$ is isogeneous to the self-product $E\times E$ for some elliptic curve $E/K$ with complex multiplication can be characterized by a property of the form
$q_C$; see Theorem \ref{th:prod} below. For convenience, let us call a curve $C/K$ with this property a
\emph{curve of CM product type}. In this case the invariant $q_C$ can be viewed as an equivalence class of positive integral \emph{ternary} quadratic forms. 

For such a curve $C/K$, it turns out that there are only finitely many isomorphism classes of curves 
$C'/K$ which have the ``same'' refined Humbert invariant. The purpose of this paper is to determine this number precisely.   

To state the result, let $d_C = d(q_C)$ denote the discriminant of $q_C$ (in the sense of Brandt \cite{B1} or of Watson \cite{Wa}), and let $\De_C = \De_{q_C}:= d_C/16$. Furthermore, let 
$\kappa_C = \kappa_{q_C} \ge1$ be defined as follows. If $q_C$ is a primitive form, then 
$\kappa_{q_C} = -I_1(q_C)/16$, where $I_1(q)$ is the genus invariant of a ternary form $q$ (as defined by Brandt \cite{B1}; see \S\ref{s:forms}), and if $q_C$ is imprimitive, then 
$\kappa_{q_C} = -I_1(q_C/4)$.                

\begin{thm}
\label{th:main} 
If $C/K$ is a curve of genus $2$ of CM product type, then the number $N_C$ of isomorphism classes of genus $2$ curves $C'/K$ whose refined Humbert invariant $q_{C'}$ is equivalent to that of $C$ is given by the formula
\beq
\label{eq:main}
N_C \;=\; 2^{\omega(\kappa_C)}h(\Delta_C)\frac{|\Aut(C)|}{|\Aut(q_C)|},  
\eeq  
where $\omega(\kappa_C)$ denotes the number of distinct prime divisors of $\kappa_C$, and $h(\De_C)$ 
denotes the number of proper equivalence classes of positive primitive binary forms of 
discriminant $\De_C$. 
\end{thm}

By using standard finiteness results, the above theorem implies the following interesting fact.

\begin{cor}
\label{c:main}
For any algebraically closed field $K$ and any integer $n\ge 1$, there are only finitely many isomorphism classes of genus $2$ curves $C/K$ of CM product type such that $N_C \le n$. Moreover, there exist infinitely many isomorphism classes of curves $C/K$ of CM product type such that $N_C > n$.  
\end{cor}

\iffalse
The cardinality of the set $\mathcal{H}(G,m)$, which consists of isomorphism classes of curves of CM product type with automorphism group $G$ admitting an elliptic subcover of degree $m$, has been studied in several cases (cf.\ Shaska \cite{Shaska_subcovers} and \cite{Ki_Aut}). Whenever $\mathcal{H}(G,m)$ is finite, its cardinality $|\mathcal{H}(G,m)|$ can be expressed as a sum of the $N_C$ associated with explicit ternary forms $q_C$ (see, e.g., Equation (6.5) of \cite{Ki_Aut}). Therefore, our result in \sref{eq:main} yields an explicit formula for $|\mathcal{H}(G,m)|$ in these cases.
\fi

It is interesting to observe that the analogue of this corollary is not true if $C/K$ is a 
\emph{superspecial curve}, i.e., if $J_C \simeq E\times E$, where $E/K$ is a supersingular elliptic curve. Such a curve can also be characterized by properties of its refined Humbert invariant $q_C$ (see Theorem \ref{th:superspecial} below). However, since by \cite{IKO} there are only finitely many isomorphism classes of superspecial curves over $K$, the number $N_C$ is bounded, contrary to the CM product case. 

The result of Theorem \ref{th:main} is a special case of a more general formula which is valid for any principally polarized abelian surfaces $(A,\theta)$ of CM product type. For this, note that the refined Humbert invariant $q_{(A,\theta)}$ is defined for any principally polarized abelian surface $(A,\theta)$ 
(cf.\ \S\ref{s:forms}),  and that we have by definition that $q_C := q_{(J_C,\theta_C)}$, if $C/K$ is a curve of genus $2$.    

In order to generalize formula \sref{eq:main} to a principally polarized abelian surface $(A,\theta)$, we first observe that all the quantities on the right hand of \sref{eq:main} naturally generalize to more general ternary forms, with the exception of the factor $|\Aut(C)|$ (cf.\ \S\ref{s:forms} below). But in \cite{K-CAS} it was shown that $|\Aut(C)|$ is given by an expression involving the number $r_n^*(q_C)$ of primitive representations of an integer $n$ by $q_C$. More precisely, we have by
Theorem 25 of \cite{K-CAS} that $|\Aut(C)| = 2a(q_C)$, 
where for any form $q$ we put
\beq
\label{eq:main2}
a(q)\; :=\; \max(1, r^*_1(q))\max(1, r^*_4(q), 3r^*_4(q) - 12). 
\eeq

This leads to the following generalization of Theorem \ref{th:main}.

\begin{thm}
\label{th:main*}
If $(A,\theta)/K$ is a principally polarized abelian surface of CM product type with refined Humbert invariant $q := q_{(A,\theta)}$, then the set $\mathcal{N}_{(A,\theta)}$ of isomorphism classes 
of principally polarized abelian surfaces $(A',\theta')/K$ such that $q_{(A',\theta')}$ is equivalent 
to $q$ is finite. Moreover, if $N_{(A,\theta)}:=|\mathcal{N}_{(A,\theta)}|$ denotes the number of these isomorphism classes, then we have that 
\beq
\label{eq:main*1}
N_{(A,\theta)} \;=\;  2^{\omega(\kappa_q) + 1}h(\Delta_q)\frac{a(q)}{|\Aut(q)|},  
\eeq
except when $q$ is equivalent to $x^2 +4\kappa(y^2+\ve yz + z^2)$, for some $\kappa > 1$ and for some 
$\ve = 0$ or $1$. In this exceptional case we have that
\beq
\label{eq:main*2}
N_{(A,\theta)} \;=\;  (2^{\omega(\kappa)-1} + 1 + \ve)\frac{h(-(4-\ve)\kappa^2)}{2+\ve}.  
\eeq 
\end{thm}   

Note that the exceptional case of Theorem \ref{th:main*} cannot occur when $(A,\theta) = (J_C,\theta_C)$ is a Jacobian; cf.\ \S3 below.    

In \S\ref{s:examples} we present some explicit examples which show that in certain cases the curve 
$C/K$ is uniquely determined by its refined Humbert invariant $q_C$; cf.\ Proposition \ref{p:example} and Examples \ref{ex:1} and \ref{ex:2}.

It is perhaps useful to mention that when char$(K) = 0$, the set $\mathcal{N}_{(A,\theta)}$ of isomorphism classes as in Theorem 3 %of principally polarized abelian surfaces 
can be identified with a \emph{generalized Humbert scheme} $H(q)$, which is attached to an arbitrary
quadratic form $q$. (These $H(q)$'s  were introduced in \S3 of \cite{K-MJ} and generalize the usual Humbert surfaces $H_N$.) Thus $\mathcal{N}_{(A,\theta)} = H(q_{(A,\theta})$, as was shown in Proposition 7.2.3 of \cite{kirPHD}. Using this fact and Theorem \ref{th:main*}, one can count 
(cf.\ \cite{kirPHD}) the number of \emph{CM-points} of a certain type and, in particular, derive a formula for $|{\cal H}(G,m)|$, where ${\cal H}(G,m)$ consists of the isomorphism classes of curves of CM product type with automorphism group $G$ admitting an elliptic subcover of degree $m$. The set 
${\cal H}(G,m)$ has been studied in several cases (cf.\ \cite{Shaska_subcovers} and \cite{Ki_Aut}).

\medskip
\noindent \textbf{Acknowledgment.} This research was partially supported by an NSERC Discovery Grant held by the first author.

%****************************
\section{The quadratic forms \protect{$q_{A}$} and \protect{$q_{(A,\theta)}$}} 
\label{s:forms}

Let $A/K$ be an abelian surface, where $K$ is an arbitrary algebraically closed field and let 
$\NS(A) = \Div(A)/\!\equiv$ denote its N\'eron-Severi group. 
(Here $\equiv$ denotes the numerical equivalence of divisors on $A$.) The intersection product 
$(D_1.D_2)$ of divisors $D_1, D_2$ on $A$ defines an integral quadratic form $q_A$ on $\NS(A)$ which is given by the formula $q_A(D) = \frac12(D.D)$, for $D \in \NS(A)$. Since $\NS(A) \simeq \Z^\rho$,
where $\rho = \rho(A)$ is the Picard number of $A$, the form $q_A$ is equivalent to an integral   quadratic form $q$ in $\rho$ variables, i.e.\ we have an isomorphism 
$(\NS(A), q_A) \simeq (\Z^\rho, q)$ of quadratic modules. We then write $q_A \sim q$. 

Now suppose that $A$ has a principal polarization $\theta\in \cP(A)$. (Here and below, $\cP(A) \subset \NS(A)$ denotes the set of principal polarization on $A$.) Then the quadratic form
$\tilde q_{(A,\theta)}$ on $\NS(A)$ is defined by the formula
\beq
\label{eq:qth} 
\tilde q_{(A,\theta)}(D) \;=\; (D.\theta)^2 - 2(D.D),\quad\mbox{for } D\in \NS(A).     
\eeq 
It is easy to see that $\tilde q_{(A,\theta)}(D+n\theta) = q_{(A,\theta)}(D)$, for all 
$n\in\Z$, so $\tilde q_{(A,\theta)}$ induces a quadratic form $q_{(A,\theta)}$ on the quotient module 
$$
\NS(A,\theta) \;=\; \NS(A)/\Z\theta.
$$  
Moreover, the Hodge Index Theorem shows that $q_{(A,\theta)}$ is a positive-definite form on $\NS(A,\theta) \simeq \Z^{\rho-1}$; cf.\ \cite{K-EC}. 
The quadratic form $q_{(A,\theta)}$ or, more correctly, the quadratic module 
$(\NS(A,\theta), q_{(A,\theta)})$ is called the \emph{refined Humbert invariant} of the principally polarized abelian surface $(A,\theta)$; cf.\ \cite{K-MJ}, \cite{K-JT}.   

We observe that if $(A,\theta) \simeq (A',\theta')$, i.e., if we have an isomorphism 
$\vp: A\stackrel\sim\ra A'$ of abelian surfaces such that $\vp^*(\theta') = \theta$, then we have an
induced module isomorphism $\vp^* : \NS(A',\theta') \stackrel\sim\ra \NS(A,\theta)$ such that 
$q_{(A',\theta')} = q_{(A,\theta)}\circ \vp^*$, and hence the associated quadratic modules are isomorphic. Equivalently, this means that $q_{(A,\theta)} \sim q_{(A',\theta')}$. 

%******
% new 21.03.26   
For later use we mention here that in the case of an abelian product surface $A = E_1\times E_2$ we have by \cite{K-MJ}, Proposition 23, an isomorphism
\beq
\label{eq:D} 
\mathbf{D} : \Z\times \Z\times \Hom(E_1,E_2)\; \stackrel\sim\ra\; \NS(E_1\times E_2) 
\eeq
of abelian groups such that for all $a_i, b_i\in\Z$ and $f_i\in \Hom(E_1,E_2)$, $i=1,2$,
%\iffalse
\beq
\label{eq:inter} 
(\mathbf{D}(a_1,b_1,f_1).\mathbf{D}(a_2,b_2,f_2)) \;=\; a_1b_2 + a_2b_1 - \be_d(f_1,f_2), 
\eeq
where $\be_d(f_1, f_2) = \deg(f_1+f_2) - \deg(f_1) - \deg(f_2)$ is the bilinear form defined by the degree map. 
%******

Suppose now that $A$ is isogeneous to a product abelian surface $E\times E$, where $E/K$ is an elliptic curve with complex multiplication, i.e. $\End(E)$ is isomorphic to an order in an imaginary quadratic field. The set of such abelian surfaces may be classified as follows.  

\begin{thm}
\label{th:prod}
If $A/K$ is an abelian surface, then the following conditions are equivalent:
\sps

\emph{(i)} $A$ is isogeneous to a product $E\times E$, for some CM elliptic curve $E/K$.
\spm

\emph{(ii)} $A$ is isomorphic to a product $E_1\times E_2$, where $E_1/K$ and $E_2/K$ are two isogeneous CM elliptic curves over $K$.
\spm

\emph{(iii)} $\cP(A) \ne \emptyset$, and for one (and hence for any) $\theta\in\cP(A)$, the refined Humbert invariant $q_{(A,\theta)}$ is a ternary form such that $q_{(A,\theta)}(D) = n^2$ is a positive square, for some $D\in\NS(A,\theta)$.
\end{thm}

\pf\ (i) $\Ra$ (ii): If $K=\C$, then this is due to Shioda and Mitani \cite{SM}. For a general ground field $K$, this is a special case of Theorem 2 of \cite{K-PCM}.  

(ii) $\Ra$ (iii): Let $\theta_{E_1,E_2}$ denote the product polarization on $A':=E_1\times E_2$, i.e., 
 $\theta_{E_1,E_2}$ is the image of the divisor $E_1\times\{0\} + \{0\}\times E_2 \in\Div(A')$ in 
$\NS(A') = \Div(A')/\!\equiv$. Then $\theta_{E_1,E_2}\in\cP(A')$, and so $\cP(A)\ne\emptyset$ because if $\vp : A \stackrel\sim\ra A'$ is any isomorphism, then 
$\vp^*(\theta_{E_1,E_2})\in\cP(A)$.  

Now let $\theta \in\cP(A)$. Then $\rank(\NS(A,\theta)) = \rho(E_1\times E_2) - 1 
= \rank(\Hom(E_1,E_2)) + 1$ by \sref{eq:D}.  %Corollary 24 of \cite{K-MJ}. 
Since $E_1\sim E_2$ are isogeneous CM elliptic curves by (ii), it follows that  $\rank(\Hom(E_1, E_2)) = 2$, so 
$\rank(\NS(A,\theta)) = 3$, and hence $q_{(A,\theta)}$ is a ternary form.  

Moreover, by Theorem 1.5 of \cite{K-EC} we know that if $E\le A$ is an elliptic subgroup on $A$, then 
$q_{(A,\theta)}([E]) = (E.\theta)^2 >0$, where $[E]\in\NS(A,\theta)$ denotes the image of $E$ in $\NS(A,\theta)$. Since $E = \vp^{-1}(E_1\times\{0\})\le A$, we see that (iii) holds.      

(iii) $\Ra$ (i): Let $\theta\in\cP(A)$, and let $D\in\NS(A,\theta)$ be such that 
$q_{(A,\theta)}(D) = n^2>0$. Then $D\ne0$, so we can write $D = mD'$, where $m\ge1$ and $D'$ is primitive in $\NS(A,\theta)$, and hence $q_{(A,\theta)}(D') = \left(\frac nm\right)^2$. 
By Theorem 1.5 of \cite{K-EC} there exists an elliptic subgroup $E \le A$ such that $[E] = D'$. Thus, 
by Poincar\'e, $A$ is isogeneous to $E \times A'$, for some abelian subvariety $A'$ of $A$ 
(cf.\ \cite{Mu}, p.\ 173). Since $\dim(A) = 2$, it follows that $\dim(A') = 1$, so $A' = E'$ is an elliptic curve, and hence $A \sim E\times E'$.           

By hypothesis, $\rank(\NS(A,\theta)) = 3$, so $4 = \rho(A) = \rho(E\times E') = 
\rank(\Hom(E, E')) + 2$, i.e., $\rank(\Hom(E, E'))=2$. This means that $E$ and $E'$ are isogeneous CM curves, so $A$ satisfies (i).

\spm\sps

In view of property (ii), we call any surface satisfying these properties a \emph{CM abelian product surface}. (In the introduction these were called abelian surfaces of CM product type.) Note that property (iii) gives a characterization of such surfaces in terms of the refined Humbert invariant 
$q_{(A,\theta)}$.     

%********
%new material (inserted 21.03.26)

While the exact analogue of Theorem \ref{th:prod} is not true for superspecial abelian surfaces, we do have the following result.

\begin{thm}
\label{th:superspecial}
If $A/K$ is an abelian surface, and if $p=\Char(K) > 0$, then the following conditions are equivalent:
\sps

\emph{(i)} $A$ is isomorphic to a product $E\times E$, for some supersingular elliptic curve \spm $E/K$.

\emph{(ii)} $A$ is isomorphic to a product $E_1\times E_2$, where $E_1/K$ and $E_2/K$ are two supersingular elliptic curves over $K$.
\spm

\emph{(iii)} $\cP(A) \ne \emptyset$, and for one (and hence for any) $\theta\in\cP(A)$, the refined Humbert invariant $q_{(A,\theta)}$ is a quadratic form of rank $5$ such that $q_{(A,\theta)}(D) = n^2$, for some $D\in\NS(A,\theta)$ and some $n\in\Z$ with $p\nmid n$.
\end{thm}

\pf\ (ii) $\Ra$ (i): This follows from Deligne's Theorem (see [IKO], p.\ 145.) 
\sps 

(i) $\Ra$ (iii): This is partially similar to the proof of (ii) $\Ra$ (iii) of Theorem \ref{th:prod}. 
Let $\theta_{E,E}$ denote the product polarization on $A':=E\times E$. Then $\theta_{E,E}\in\cP(A')$, and so $\cP(A)\ne\emptyset$ because we have an isomorphism $\vp : A \stackrel\sim\ra A'$.

Now let $\theta \in \cP(A)$ be any principal polarization. Then $\rank(\NS(A,\theta)) = 
\rho(E\times E) - 1 = \rank(\End(E)) + 1 = 5$ by \sref{eq:D} and by the fact that $E$ is supersingular.

To verify the last property of (iii), we will use Corollary 24 of \cite{K-MJ} which states that 
$\theta = \mathbf{D}(a_1,a_2,f)$, for some $a_1,a_2 \in\Z$ and $f\in \End(E)$ such that 
$a_1 > 0$ and $a_1a_2 - d = 1$, where $d = \deg(f)$. We now consider two cases:
\sps

\noindent{\bf Case 1:} $p|\gcd(a_1,a_2)$.     
\sps

\noindent Then $p\nmid d$. Consider $D = \mathbf{D}(1,d,f)$. By \sref{eq:inter} we have have that
$(D.D) = 0$ and $n:=(D.\theta) = a_1d + a_2 - 2d \equiv -2d \not\equiv 0 \Mod{p}$ because $p \ne 2$ 
by \cite{IKO}, Theorem 3.3 (III). Thus by \sref{eq:qth} we have that $q_{(A,\theta)}(D) = n^2$ and $p\nmid n$, so (iii) holds in this case.
\sps

\noindent{\bf Case 2:} $p\nmid\gcd(a_1,a_2)$.     
\sps

\noindent In this case either $p\nmid a_1$ or $p\nmid a_2$. Then put $D_1 = \mathbf{D}(0, 1, 0)$, 
if $p\nmid a_1$ and put $D_2 = \mathbf{D}(1,0,0)$, if $p\nmid a_2$. Then by \sref{eq:inter} we have that
$(D_i.D_i) =0$, and $n_i:=(D_i.\theta) = a_i \not\equiv 0 \Mod p$. Thus 
$q_{(A,\theta)}(D_i) = n_i^2$, and so (iii) holds in this case as well.        
\sps

(iii) $\Ra$ (ii):  Let $\theta\in\cP(A)$, and let $D\in\NS(A,\theta)$ be such that 
$q_{(A,\theta)}(D) = n^2 \not\equiv0\Mod p$. Then $D\ne0$, so we can write $D = mD'$, where $m\ge1$ and $D'$ is primitive in $\NS(A,\theta)$, and hence $q_{(A,\theta)}(D') = N^2$ with 
$N:=\left(\frac nm\right)^2$. 
By Theorem 18 of \cite{K-ESC1} there exists an $N$-presentation $(E_1, E_2, \psi, \pi)$ of $A/K$. This means that we have an isogeny $\pi : E_1 \times E_2\ra A$ and an anti-isometry $\psi: E_1[N] \ra E_2[N]$
such that $\Ker(\pi) = \mbox{Graph}(-\psi)$. 

By hypothesis, $\rho(A) = \rank(\NS(A,\theta)) + 1 = 6$. Since $\pi$ is an isogeny, we have that
$\rho(A) = \rho(E_1\times E_2) = \rank(\Hom(E_1,E_2)) + 2$ by \sref{eq:D}, and so 
$\rank(\Hom(E_1,E_2)) = 4$. This means that $E_1$ and $E_2$ are supersingular elliptic curves. 

Since $p\nmid N$, we have that $|\Hom(E_1[N],E_2[N])| = N^4$. Moreover, since 
${\cal H}:= \Hom(E_1, E_2) \simeq \Z^4$, we also have that $|{\cal H}/N{\cal H}| = N^4$, and so 
Lemma 22 of \cite{K-ESC1} shows that the map
$\rho_N: \Hom(E_1,E_2) \ra \Hom(E_1[N],E_2[N])$ (given by $\rho_N(h) = h_{|E_1[N]}$) is surjective. Thus $\psi = h_{|E[N]}$, for some $h \in \Hom(E_1[N],E_2[N])$. This means that the isogeny defect $m_\psi = 1$, and hence $A \simeq E_1 \times E_2$ by Theorem 4 
of \cite{K-ESC1}. (Note that the proof of Theorem 4 is valid for all principally polarized abelian surfaces $(A,\theta)$ and not just Jacobians.) This shows that (ii) holds.    

\spm\sps

We now come back to the case of a CM product surface, and 
%********
recall from \cite{K-MJ}, \cite{K-JT} and \cite{Ki} some basic arithmetic facts about the forms 
$q_A$ and $q_{(A,\theta)}$ for such a surface. To state these, we will use discriminant 
$\disc(f) = d(f)$ of an integral quadratic form $f$ as defined in Watson \cite{Wa}, p.~2. Moreover, the \emph{content} $\cont(f)$ of such a form is the gcd of its values. (It is easy to see that when $f$ is
a ternary form, then $\cont(f) = t$ is the ``coefficient-divisor'' of \cite{B1}.) 
If $\cont(f)=1$, then $f$ is said to be \emph{primitive}; cf.\ \cite{Wa}, p.~4. In addition, the basic invariants $I_k(f)$ for $k=1,2$ of an integral primitive ternary form $f$ are defined as in Brandt \cite{B1}. (These invariants are closely related to the invariants $\Om$ and $\De$ as defined in Dickson \cite{Di}.)

More precisely,  the genus invariant $|I_1(q)|$ of a primitive form $q$ is the content of the \emph{adjoint form} adj$(q)$, which is defined by the formula (as in  Watson \cite{Wa}, p.~ 25)
\begin{equation}
\label{eq: adj_form}
    M(\text{adj}(q)) \;=\; -2\text{adj}(M(q)) \;=\; -2\det(M(q))M(q)^{-1}.
\end{equation}
The sign is assigned so that the \emph{reciprocal} 
$F_f := \mbox{adj}(f)/I_1(f)$ of $f$ is a positive definite form when $f$ is positive definite. 
Furthermore, $I_2(f) := I_1(F_f)$.

\begin{prop} 
\label{p:invars}
Let $A \simeq E\times E'$ be a CM abelian product surface, and let $q_{E,E'}$ denote the degree form 
on $\Hom(E,E')$. Let $\De = d(q_{E,E'})$ denote its discriminant and $\kappa = \cont(q_{E,E'})$ its content. If $\theta \in \NS(A)$ is a principal polarization, then 
\beq
\label{eq:invars1}
d(q_A) \;=\; \De \quad\mbox{and}\quad d(q_{(A,\theta)}) \;=\; 16\De.    
\eeq     
Moreover, if $c:=\cont(q_{(A,\theta)})$, then $c=1$ or $4$, and we have that 
\beq
\label{eq:invars2}
I_1(q_{(A,\theta)}/c) \;=\; -\frac{16\kappa}{c^2}\quad\mbox{and}\quad
I_2(q_{(A,\theta)}/c) \;=\; \frac{c\De}{\kappa^2}.    
\eeq 
\end{prop}

\pf\ We have that $\rank(\Hom(E,E')) = 2$ because $E$ and $E'$ are isogeneous CM elliptic curves, and
so $q_{E,E'}$ is a binary quadratic form. It thus follows from Corollary 24 of \cite{K-MJ} that  
$\det(q_A) = -\det(q_{E,E'})$, where (as in \cite{K-MJ}) $\det(f)$ denotes the determinant of the associated Gram matrix $M(f)$. Thus, since $\rank(\NS(A)) =4$, it follows from \cite{Wa}, p.\ 2, that  
$d(q_A) = \det(q_A) = - \det(q_{E,E'}) = d(q_{E,E'})$. This proves the first equality of 
\sref{eq:invars1}.         

The second equality of \sref{eq:invars1} follows from Proposition 9 of \cite{K-MJ}. Indeed, by that result we have that $\det(q_{(A,\theta)}) = -2^5\det(q_A)$, so 
$d(q_{(A,\theta)}) = -\frac12\det(q_{(A,\theta)}) = 2^4\det(q_A) = 16d(q_A) = 16\De$. This proves
 \sref{eq:invars1}. 

If $c=1$, i.e., if $q_{(A,\theta)}$ is primitive, then formula \sref{eq:invars2} follows from
Propositions 25 and 18 of \cite{K-rH}. If $c >1$, then 
by Corollary 3.5 of \cite{Ki} and the proof of Proposition 4.1 of \cite{Ki}, we see that $c=4$, and so 
in this case \sref{eq:invars2} follows directly from \spm\sps
that proposition.   

In view of the above result, it is useful to introduce the following notation (which was already used in 
Theorem \ref{th:main*}).
\spm\sps

\notn\ If $q$ is a positive ternary quadratic form such that $c:=\cont(q) = 1$ or $4$, then put 
$$
\De_q\;:=\; d(q)/16\;\quad\mbox{and}\quad \kappa_q \;:=\;-c^2I_1(q/c)/16. 
$$

\spm\sps

It is clear from the definition that the form $q_{(A,\theta)}$ is determined by the form $q_A$ and 
an element $\theta$ with $q_A(\theta) = 1$; this is the so-called $\theta$-construction in 
\cite{K-JT}. We now show conversely that the ternary form $q_{(A,\theta)}$ determines the quaternary 
form $q_A$ (up to equivalence). More precisely, we prove the following result.

\begin{thm}
\label{th:qA}
Let $A_1$ and $A_2$ be two CM abelian product surfaces, and let $\theta_i\in \NS(A_i)$ be a principal polarization on $A_i$, for $i=1,2$. If $q_{(A_1,\theta_1)}$ is equivalent to $q_{(A_2,\theta_2)}$, 
then $q_{A_1}$ is equivalent to $q_{A_2}$.    
\end{thm} 

To prove this result, we will use several facts which were proven elsewhere. The first concerns the 
$p$-adic equivalence class of the form $q_{(A,\theta)}$, where the $p$-adic equivalence of two integral forms $f_1$ and $f_2$ is defined as in Jones \cite{Jo}, p.\ 82, and is denoted by $f_1 \sim_p f_2$. Recall also that $f_1$ and $f_2$ are \emph{genus-equivalent} if $f_1 \sim_p f_2$, for all primes $p$ (including $p=\infty$). We let $\gen(f_1)$ denote the set of (equivalence classes of) integral forms 
$f_2$ which are genus-equivalent to $f_1$. 
    
\begin{prop}
\label{p:p-adic}
Let $A = E\times E'$ be a CM abelian product surface, and let $\theta$ be a principal polarization on 
$A$. If $f_q := x^2 \perp 4q$, where $q=q_{E,E'}$ is the degree form, then 
\beq
\label{eq:p-adic1}
q_{(A,\theta)} \;\sim_p\; f_q,\quad \mbox{for all odd primes $p$.} 
\eeq
Furthermore, if $q_{(A,\theta)}$ is primitive, then \emph{\sref{eq:p-adic1}} also holds for $p=2$, and so in this case $q_{(A,\theta)}$ is genus-equivalent to $f_q$.
\end{prop}
 
\pf\ The first assertion \sref{eq:p-adic1} follows immediately from Corollary 19 of \cite{K-JT} because by equation (29) of \cite{K-rH} we know that $f_q \sim q_{(A,\theta_{E,E'})}$, 
where $\theta_{E,E'}$ is the product polarization on $A$. The second assertion follows from Theorem 20 of \cite{K-JT}, as is explained in the proof of Proposition 25 of \cite{K-rH}.   
\spm\sps

In order to state the next result, we recall from \cite{K-rH} and \cite{Ki} that if $f: M \ra \Z$ is a quadratic form on a module $M$, then $R(f) := \{f(x): x\in M\}$ denotes the set of values represented by $f$, and if $a, m$ are integers, then 
$$
R_{a,m}(f) \;:=\;  \{n\in R(f): n\equiv a \Mod m\}.
$$  
 
\begin{prop}
\label{p:imp}
Let $A = E\times E'$ be a CM abelian product surface, and let $\theta$ be a principal polarization on 
$A$. If $q_{(A,\theta)}$ is imprimitive, then $d(q_{E,E'}) \equiv 0\Mod 4$ and 
$R_{3,4}(q_{E,E'}) \ne \emptyset$. In particular, $\kappa := \cont(q_{E,E'})$ is odd. Furthermore, there exists $n \in R_{3,4}(q_{E,E'})$ such that $\frac{n}{\kappa} \in R(F_{q_{(A,\theta)}/4})$.   
\end{prop}

\pf\ The first assertions follow immediately from Proposition 3.6 of \cite{Ki}. Moreover, the last
assertion follows from Propositions 3.7 and 4.1 of \cite{Ki}. To see this, put $n = q_{E, E'}(h)$, where
$h\in\Hom(E,E')$ is as in Proposition 4.1 of \cite{Ki}. 
Then by Propositions 3.7 and 4.1 of \cite{Ki} we have that 
$n \in R_{3,4}(q_{E,E'})$ and that $\frac{n}{\kappa} \in R(F_f)$, where $F_f$ denotes the reciprocal of the primitive form $f := q_{(A,\theta)}/4$.    
\spm\sps

We are now ready to prove Theorem \ref{th:qA}.
\spm\sps

\pff\ Theorem \ref{th:qA}. For $i=1,2$, let $E_i/K$ and $E'_i/K$ be two (isogeneous) CM elliptic curves such that $A_i \simeq E_i\times E'_i$. We first observe that it suffices to show that the given hypothesis implies that the binary forms $q_1:=q_{E_1,E'_1}$ and $q_2:=q_{E_2,E'_2}$ are genus-equivalent, i.e., that 
\beq
\label{eq:qA1} 
\gen(q_1) \;=\; \gen(q_2).  
\eeq
Indeed, since we have that
\beq
\label{eq:q_A2}
q_{A_i} \;\sim\; xy \perp (-q_i), \quad\mbox{for } i=1,2, 
\eeq
by \cite{K-MJ}, Proposition 23 (or by  \cite{K-JT}, formula (6)),
it follows from Remark 27 of \cite{K-JT} that \sref{eq:qA1} implies that 
$q_{A_1} \sim q_{A_2}$, as desired.  

To prove that \sref{eq:qA1} holds, we will distinguish two cases.
\spm

\noindent{\bf Case 1.} $q:=q_{(A_1,\theta_1)}$ is primitive.   
\sps

\noindent Since $q$ is primitive, so is $q_{(A_2,\theta_2)} \sim q$. Thus, by Proposition \ref{p:p-adic}
we have that $q_{(A_i,\theta_i)} \in \gen(f_{q_i})$, for $i=1,2$.
Thus, $\gen(f_{q_1}) = \gen(f_{q_2})$ because $q_{(A_1,\theta_1)} \sim q_{(A_2,\theta_2)}$, and so
$x^2 + 4q_1 \sim_p x^2 + 4q_2$, for all primes $p\ge2$. Since $x^2$ is a form of unit determinant 
in the sense of Jones \cite{Jo}, it follows from the Cancellation Theorem 37 of \cite{Jo}
that $2q_1 \sim_p 2q_2$, for all primes $p\ge2$, and so $\gen(2q_1) = \gen(2q_2)$ because $2q_1$ and 
$2q_2$ are both positive binary forms. Then also $\gen(q_1) = \gen(q_2)$, which proves \sref{eq:qA1} in this case.    
\spm

\noindent{\bf Case 2.} $q =q_{(A_1,\theta_1)}$ is imprimitive.   
\sps

\noindent Since $q$ is imprimitive, so is $q_{(A_2,\theta_2)} \sim q$, and hence both have content 
$c=4$ by Proposition \ref{p:invars}. Since $I_k(q_{(A_1,\theta_1)}/4) = I_k(q_{(A_2,\theta_2)}/4)$, 
for $k=1,2$, it therefore follows from \sref{eq:invars2} that $q_1$ and $q_2$ have the same content
$\kappa$ and the same discriminant $\De$. Thus, $q'_i := q_i/\kappa$ are two primitive binary forms 
of the same discriminant $\De' = \De/\kappa^2$. It is clear that \sref{eq:qA1} follows once we have shown that $\gen(q'_1) = \gen(q'_2)$.

To prove this, we will use Theorem 3.21 and Lemma 3.20 of Cox \cite{Cox} which state that 
$\gen(q'_1) = \gen(q'_2)$ if and only if $q'_1$ and $q'_2$ have the same ``complete character''. 

This latter condition is defined as follows. For an odd prime $p$, let $\chi_p$ denote the Legendre character which is defined by $\chi_p(a) = \leg ap$, for $a\in\Z$, and put $\chi_{-4}(n) = \delta(n) = (-1)^{(n-1)/2}$ and $\chi_{8}(n) = \epsilon(n) = (-1)^{(n^2-1)/8}$, for $n\equiv 1\Mod 2$. 
(Here and below we will use the $\chi_n$ notation of \cite{B1} instead of the $\de,\epsilon$ notation of \cite{Cox}.) 
Then the list of \emph{assigned characters} of discriminant $\De'$ is 
$X(\De') := \{\chi_p : p|\De'\} \,\cup\, X_s(\De')$, where 
$X_s(\De') \subset X_s:= \{\chi_{-4}, \chi_8, \chi_{-4}\chi_8\}$ is the set of \emph{supplementary characters} which is given in the table on p.\ 55 of \cite{Cox}. For each $\chi \in X(\De')$, 
the value $\chi(q'_i) := \chi(r_i)$ does not depend on the choice of $r_i\in R(q'_i)$, provided that
$\gcd(r_i,\De') = 1$. Thus, by Cox \cite{Cox}, loc.\ cit.,  we have that $\gen(q'_1) = \gen(q'_2)$ if and only if          
\beq
\label{eq:qA3} 
\chi(q'_1)\;=\; \chi(q'_2),\quad\mbox{for all } \chi \in X(\De').  
\eeq

It thus suffices to verify \sref{eq:qA3}. For this, we recall that there is a similar theory of assigned characters for primitive ternary forms; cf.\ Smith \cite{S} and Brandt \cite{B1}, \cite{B2}.

We will first apply this theory to the ternary forms $f_{q_i}$, for $i=1,2$.  Since $f_{q_i}$ is primitive and $f_{q_i} \sim q_{(A_i, \theta_{E_i,E'_i})}$ (cf.\ the proof of 
Proposition \ref{p:p-adic}), it follows from Proposition \ref{p:invars} that
\beq
\label{eq:qA4}
d(f_{q_i}) \;=\; 16\De,\quad I_1(f_{q_i}) = -16\kappa\quad\mbox{and}\quad
I_2(f_{q_i}) = \De', \quad\mbox{for $i=1,2$.}    
\eeq

Let $p$ be an odd prime. Since $q_{(A_1,\theta_1)} \sim q_{(A_2,\theta_2)}$ by hypothesis, 
it follows from \sref{eq:p-adic1} that $f_{q_1} \sim_p q_{(A_1,\theta_1)} \sim_p q_{(A_2,\theta_2)} \sim_p f_{q_2}$. This (together with \sref{eq:qA4}) implies that 
their reciprocals $F_i := F_{f_{q_i}}$ are also $p$-adically equivalent (for $p$ odd), as is well known; 
cf.\ Proposition 3.2.10(iii) of \cite{kirPHD}. Furthermore, from that proposition we also obtain that 
\beq
\label{eq:qA5}
\chi_p(n_1) = \chi_p(n_2),\; \mbox{$\forall$ odd $p|\De'$, $\forall n_i \in R(F_i)$ with 
$\gcd(n_i,\De') =1$, $i=1,2$.}   
\eeq
Now by formula (44) of \cite{K-JT} we have that $F_i \sim (-\De'\kappa)x^2 \perp q'_i$, so 
$R(q'_i) \subset R(F_i)$, for $i=1, 2$. Thus, if $p$ is an odd prime with $p|\De'$, 
then $\chi_p\in X(\De')$, so if $n_i\in R(q'_i)$ with $\gcd(n_i,\De') = 1$, then 
by \sref{eq:qA5} we obtain that 
$\chi_p(q'_1) =\chi_p(n_1) = \chi_p(n_2) = \chi_p(q'_2)$. This means \sref{eq:qA3} holds for all
$\chi \in X(\De') \setminus X_s(\De')$. 

It thus remains to show that  \sref{eq:qA3} also holds for $\chi \in  X_s(\De')$. For this, note first that by Proposition \ref{p:imp} we have that $\De' \equiv \De \equiv 0\Mod4$. As in \cite{Cox}, p.\ 55, put $n:= -\frac{\De'}{4}$. 

Suppose first that $n\equiv3\Mod4$. Then $X_s(\De') = \emptyset$ by \cite{Cox}, so there is nothing to prove. Next, suppose that $n\equiv1\Mod4$, so $X_s(\De') = \{\chi_{-4}\}$ by \cite{Cox}. Now by Proposition
\ref{p:imp} there exists $n_i\in R_{3,4}(q_i)$, so $n'_i := \frac{n_i}\kappa \in R_{1,2}(q_i)$ because
$\kappa$ is odd. Thus, for $i=1,2$, we have that
$\chi_{-4}(q'_i) = \chi_{-4}(n_i') = \chi_{-4}(\kappa)^{-1}\chi_{-4}(n_i) = -\chi_{-4}(\kappa)$ because $n_i\equiv 3\Mod4$, and so $\chi_{-4}(q'_1) = \chi_{-4}(q'_2)$. Thus, \sref{eq:qA3} holds in this case as well.

We are thus left with the case that $n\equiv0\Mod2$, i.e., that $\De'\equiv 0\Mod 8$. In this case we will make use of the supplementary characters of the reciprocals $F_{f_i}$ of the forms 
$f_i = q_{(A_i,\theta_i)}/4$, for $i=1,2$. By hypothesis, $f_1 \sim f_2$, and by 
Proposition \ref{p:invars} we see that $f_i$ is a primitive ternary form with genus invariants 
$I_1(f_i) = -\kappa$ and $I_2(f_i) = 4\De'$. 
Since $\De' \equiv 0 \Mod8$, we have that $I_1(F_{f_i}) = I_2(f_i) \equiv 0\Mod{32}$. 
By Brandt \cite{B1}, p.\ 337, this implies that $\chi_{-4}$, $\chi_{8}$ and 
$\chi_{-8} := \chi_{-4}\chi_8$ are assigned characters of $F_{f_i}$. This means that 
if $\chi \in X_s = \{\chi_{-4}, \chi_{8}, \chi_{-8}\}$, and if $i=1,2$, then we have that
\beq
\label{eq:qA6}
\chi(r) \;=\; \chi(r'), \quad\mbox{for all } r,r' \in R_{1,2}(F_{f_i}).
\eeq    

Now let $\chi \in X_s(\De') \subset X_s$. Then $\chi(q'_i) = \chi(n_i)$,
for any $n_i \in R_{1,2}(q'_i)$. By Proposition \ref{p:imp} there exists 
$n_i \in  R_{1,2}(q'_i) \cap R_{1,2}(F_{f_i})$. Now since $f_1\sim f_2$, we also have that 
$F_{f_1}\sim F_{f_2}$ and so $R_{1,2}(F_{f_1})= R_{1,2}(F_{f_2})$, and so 
$n_1,n_2\in R_{1,2}(F_{f_1})$. Thus, by \sref{eq:qA6} we obtain that
$\chi(q'_1) = \chi(n_1) = \chi(n_2) = \chi(q'_2)$, for all $\chi \in X_s(\De)$. This proves that
\sref{eq:qA3} holds in all cases, and hence also \sref{eq:qA1}, as desired.

%**********************

\section{The number \protect{$N_{(A,\theta)}$} of isomorphism classes}
\label{s:formulas}

The purpose of this section is to prove the main results stated in the introduction. Let $\langle A,\theta\rangle$ denote the isomorphism class of a principally polarized abelian surface $(A,\theta)$, 
and let $\mathcal{A}_2(K)$ denote the set of isomorphism classes of principally polarized 
abelian surfaces $(A,\theta)/K$. When $(A,\theta)/K$ is a fixed principally polarized abelian surface, %where $A$ is a CM abelian product surface, let us 
put
$$
\mathcal{N}_{(A,\theta)}\;:=\;\{\langle A', \theta' \rangle\in \mathcal{A}_2(K)\, :\,  q_{(A',\theta')}\sim q_{(A,\theta)}\}.
$$
By Theorem \ref{th:prod}, it is clear that if $A$ is a CM abelian product surface, then so is $A'$, 
for any $\langle A',\theta'\rangle \in \mathcal{N}_{(A,\theta)}$. As in the introduction, let us put
$$
 N_{(A,\theta)} \;=\; |\mathcal{N}_{(A,\theta)}|.
$$

The first step is to count the number of non-isomorphic abelian surfaces $A'/K$, which appear in the set $\mathcal{N}_{(A,\theta)}$. To formulate this, let $\langle A \rangle$ denote the 
isomorphism class of an abelian surface $A/K$.  For a given quadratic form $q$, let us define 
$$
\mathcal{N}(q)\;:=\;\{\langle A \rangle \, :\,  q_{(A,\theta)}\sim q, 
\mbox{ for some }\theta\in\mathcal{P}(A)\}
$$
as the set of isomorphism classes of abelian surfaces $A/K$ such that $q_{(A,\theta)}$ is equivalent 
to $q$, for some $\theta\in\mathcal{P}(A)$. Note that if $q_{(A,\theta)}$ is a ternary form which represents a square, then we see that  $\mathcal{N}(q_{(A,\theta)})$ is a subset of the set of isomorphism classes of CM abelian product surfaces by Theorem \ref{th:prod}. 

In order to find $|\mathcal{N}(q_{(A,\theta)})|$ in this case, recall from the previous section that the form $q_{(A,\theta)}$ determines the
form $q_A$ (up to equivalence);  cf.\ Theorem \ref{th:qA}. This fact is the key tool to complete the first step because it provides the relation between the set  $\mathcal{N}(q_{(A,\theta)})$ and the set $\mathcal{N}({A})$, where 
$$
\mathcal{N}({A})\;:=\;\{ \langle A' \rangle \, :\, A' \sim A,\ q_{A'}\sim q_A\}
$$
denotes the set of isomorphism classes of abelian surfaces $A'/K$ which are isogeneous to a given $A$ and whose intersection form $q_{A'}$ is equivalent to $q_A$. The advantage of studying this set is 
that \cite{K-PCM} gives an explicit formula for the cardinality $|\mathcal{N}({A})|$ when $A$ is a 
CM abelian product surface.  %It turns out that the first step is done as follows.
We now prove the following result, which will complete the first step. 

\begin{prop}
\label{p: |A(q)|}
 Let  $A= E_1\times E_2$ be a CM abelian product surface, and let $\theta\in\mathcal{P}(A)$. Put  $q:=q_{(A,\theta)}$, and  let $q_1:=q_{E_1,E_2}$. Then we have that 
\beq
\label{eq:Nq=NA}
\mathcal{N}(q)\;=\;\mathcal{N}(A).
\eeq
Moreover, if we put $ \Delta' = d(q_1)/\cont(q_1)^2$, then  
\beq
\label{e: |A(q)|}
    |\mathcal{N}(q)|\;=\;\frac{h(\Delta')}{g(\Delta')},
    \eeq
    where $g(\Delta')$ denotes the number of genera of discriminant $\Delta'$.
\end{prop}

\pf\ We first observe that 
\beq
\label{eq:Nq=NA1}
\End^0(E_i) \;\simeq\; \Q(\sqrt{d(q_1)}),\quad\mbox{for } i=1,2. 
\eeq
Indeed, since $E_1$ and $E_2$ are two isogeneous CM elliptic curves, we have that 
$F:= \End^0(E_1) \simeq \End^0(E_2)$ is an imaginary quadratic field, and by Corollary 42 
of \cite{K-PCM} we have that $d(q_1) = \lcm(f_{E_1}, f_{E_2})^2\Delta_F$, where $\Delta_F$ is the
discriminant of $F$ and $f_{E_i}$ is the \emph{endomorphism conductor} of $E_i$ as defined in
\cite{K-PCM}. From this, \sref{eq:Nq=NA1} follows because $F \simeq \Q(\sqrt{\De_F})$.     

To prove \sref{eq:Nq=NA}, suppose first that $\langle A' \rangle \in \mathcal{N}(q)$, so by 
definition there is a $\theta' \in \mathcal{P}(A')$ such that $q_{(A',\theta')}\sim q$. Since $A$ is a CM abelian product surface, so is $A'$ by Theorem \ref{th:prod}, and hence we have by 
Theorem \ref{th:qA} that $q_A\sim q_{A'}$. Thus, in order to show that 
$\langle A' \rangle \in \mathcal{N}(A)$, it suffices to prove that $A\sim A'$.

Now since $A'$ is a CM product surface, we have that $A'\simeq E_1' \times E_2'$, 
for some isogeneous CM elliptic curves $E'_1/K$ and $E'_2/K$. If we put $q_1' := q_{E_1',E_2'}$,   
then by \sref{eq:invars1} we see that $d(q_1) = d(q)/16 = d(q_{(A',\theta')})/16 = d(q_1')$. 
Thus, equation \sref{eq:Nq=NA1} (applied to $E'_i$) shows that 
$\End^{0}(E_i')\simeq\Q(\sqrt{d(q'_1)}) = \Q(\sqrt{d(q_1)})\simeq \End^{0}(E_i)$, for $i=1,2$, 
and so it follows that $E_i\sim E_i'$; cf.\ Proposition 36 of \cite{K-PCM}. It thus follows that 
$A' \simeq E_1' \times E_2' \sim E_1 \times E_2 = A$, and hence
$\langle A' \rangle \in \mathcal{N}(A)$.

To prove the opposite inclusion, let $\langle A' \rangle \in \mathcal{N}(A)$, so we have that $q_{A'}\sim q_A$. Hence, since $\theta\in\mathcal{P}(A)$, it follows from Proposition 29 of \cite{K-rH} that $ q\sim q_{(A',\theta')}$, for some $\theta'\in\mathcal{P}(A')$, which proves that $\langle A' \rangle \in \mathcal{N}(q)$, and so \sref{eq:Nq=NA}  follows.

To prove \sref{e: |A(q)|}, we will use Corollary 70 of \cite{K-PCM}. To apply it, note first that 
there exists an elliptic curve $E/K$ with $E\sim E_i$ such that the discriminant $\De_E$ of the order $\End(E)$ equals $d(q_1)$. Indeed, by Proposition 29 of \cite{K-PCM} there exists an $E\sim E_i$ such that
$f_E = \lcm(f_{E_1},f_{E_2})$, and so by Corollary 42 of \cite{K-PCM} we have that
$d(q_1) = \lcm(f_{E_1},f_{E_2})^2\De_F$, where $F = \End^0(E)\simeq \End^0(E_i)$, and so 
$d(q_1) = f_E^2\De_F = \De_E$.   

We thus have that $A \simeq E_1\times E_2 \sim E\times E$. Moreover, by \sref{eq:q_A2} we know that 
$q_A\sim xy \perp (-q_1)$, so it follows that 
$\mathcal{N}(A) = \mathcal{N}_{q_1}:= \{A'\sim E\times E : q_{A'} \sim xy \perp (-q_1) \}/\!\simeq$. 
We thus see that formula \sref{e: |A(q)|} follows directly from that given in Corollary 70 
\spm\sps
of \cite{K-PCM}.
%\spm\sps

The second step in finding a formula for $N_{(A,\theta)}$ is to determine the number of isomorphism classes of principal polarizations $\theta'$ on a given abelian surface $A'$ with 
$\langle A'\rangle \in \mathcal{N}(q_{(A,\theta)})$. 
\iffalse
such that $q_{(A',\theta')}\sim q_{(A,\theta)}$. This is itself very delicate study topic, and it has been already worked in detail for more general cases in \cite{K-CAS}, \cite{K-PP}. 
To discuss it, let us recall from these articles that for a given quadratic form $q$ and an abelian surface $A$, the set 
\fi
For this, we will use the results of  \cite{K-CAS} and \cite{K-PP}. As in those papers, let
$$
\mathcal{P}(A,q):=\{\theta\in \mathcal{P}(A) \, :\, q_{(A,\theta)}\sim q \}
$$
denote the set of principal polarizations $\theta$ on $A$ such that $q_{(A,\theta)}$ is equivalent 
to a given quadratic form $q$. By Theorem 9 of \cite{K-CAS} we know that the automorphism group $\Aut(A)$ of $A$ acts on 
$\mathcal{P}(A,q)$, so we can consider the set  
$\ov{\mathcal{P}}(A,q):=\Aut(A)\backslash \mathcal{P}(A,q)$ of orbits of $\mathcal{P}(A,q)$ under 
this action. 
We now show that the number $N_{(A,\theta)}$ can be expressed in terms of a sum of the cardinalities of 
suitable sets $\ov{\mathcal{P}}(A_i,q)$.   

\begin{prop}
\label{p:NAth}
Let $(A,\theta)$ be a principally polarized abelian surface, where $A$ is a CM abelian product 
surface, and let $q = q_{(A,\theta)}$. If $A_1,\ldots, A_n$ is a system of representatives of 
the finite set ${\cal N}(q)$, then
\beq
\label{e: H(q)=union}
{\cal N}_{(A,\theta)} \;=\ \{\langle A_i,\theta_i\rangle : \theta_i \in{\cal  P}(A_i,q), 
\mbox{ for } 1\le i\le n\},  
\eeq
and hence 
\beq 
\label{e: |H(q)|=summation}
   N_{(A,\theta)} \;=\; \sum_{i=1}^n |\overline{\cal P}(A_i,q)|.
\eeq 
In particular, $N_{(A,\theta)} < \infty$. 
\end{prop}

\pf\ If $\langle A', \theta'\rangle\in   \mathcal{N}_{(A,\theta)}$, then by definition 
$\langle A' \rangle \in \mathcal{N}(q)$ and $\theta'\in\mathcal{P}(A',q)$ because 
$q_{(A',\theta')}\sim q_{(A,\theta)} = q$. Thus, $A'\simeq A_i$, for some $1\leq i \leq n$, 
and hence the left hand side of \sref{e: H(q)=union} is contained in the right hand side.  

To prove the opposite inclusion, let $\langle A_i, \theta_i \rangle$ be an element contained in the right hand side of \sref{e: H(q)=union}, so $\langle A_i \rangle \in \mathcal{N}(q)$ and $\theta_i \in \mathcal{P}(A_i,q)$. This implies that $q_{(A_i,\theta_i)}\sim q = q_{(A,\theta)}$, and hence 
$\langle A_i, \theta_i\rangle\in \mathcal{N}_{(A,\theta)}$, which verifies \sref{e: H(q)=union}.

To verify \sref{e: |H(q)|=summation}, recall first from Proposition \ref{p: |A(q)|} that ${\cal N}(q)$ is a finite set, so the sum on the right hand side of \sref{e: |H(q)|=summation} is a finite sum.

Next we observe that if $\langle A_i, \theta_i\rangle = \langle A_j, \theta_j \rangle$, then in particular $A_i \simeq A_j$ and so $i=j$. Thus, $\lb A_i, \theta_i\rb
= \lb A_j, \theta_j\rb$ if and only if $i=j$ and $\theta_i$ and $\theta_j$ lie in the same 
$\Aut(A_i)$-orbit of $\cP(A_i,q)$. This means that for each $i$, the number of 
$\lb A',\theta'\rb\in \mathcal{N}_{(A,\theta)}$ with $A'\simeq A_i$ is equal to 
$|\ov{\mathcal{P}}(A_i, q)|$, and so \sref{e: |H(q)|=summation} follows from \sref{e: H(q)=union}.

Since $\ov\cP(A_i, q)$ is always a finite set (cf.\ Theorem 1 of \cite{K-CAS}, together with formula (7) of \cite{K-CAS}), it follows from
\sref{e: |H(q)|=summation} that $N_{(A,\theta)}=|{\cal N}_{(A,\theta)}|$ is also finite.

\spm\sps

To analyze $N_{(A,\theta)}$ further, we
will make use of a formula for the cardinality of the set $\ov{\mathcal{P}}(A, q)$ found in 
\cite{K-CAS} and \cite{K-PP}. One of the ingredients of this formula is the quantity $a(q)$
which was mentioned in \sref{eq:main2}, and which used certain ``representation numbers'' $r_n^*(q)$. 
These numbers are defined as follows. If $n>0$ is an integer, and if $q$ is a positive ternary form,
then the \emph{number of primitive representations} of $n$ by $q$ is defined by 
$$
r_n^*(q):=|\{(x_1,x_2,x_3)\in R_n(q)\, :\, \gcd(x_1,x_2,x_3)=1\}|,
$$
where $R_n(q) := \{(x_1,x_2,x_3)\in \Z^3\, :\, q(x_1,x_2,x_3)=n \}$
denotes the set of all representations of an integer $n$ by  $q$. 

\begin{prop}
\label{p: |P(A,q)|}
Let $A= E \times E'$ be a CM abelian product surface, and suppose that $\mathcal{P}(A, q)$ is  nonempty. Let $\Delta=d(q_{E,E'})$ denote the discriminant of the degree map $q_{E,E'}$, 
and let $\kappa=\cont(q_{E,E'})$ be its content, and  put $\Delta'=\Delta/\kappa^2$. If $q$ is 
not equivalent to $x^2 +4\kappa(y^2+\ve yz + z^2)$, for any $\kappa > 1$ and $\ve = 0$ or $1$, 
then  we have that 
\beq
\label{e: |P(A,q)|}
 |\ov{\mathcal{P}}(A, q)|
\;=\;\frac{2^{\omega(\kappa)+1}g(\Delta')h(\Delta)a(q)}{|\Aut(q)|h(\Delta')}.
 \eeq
On the other hand, if $q$ is  equivalent to $x^2 +4\kappa(y^2+\ve yz + z^2)$, for some $\kappa > 1$ and for some $\ve = 0$ or $1$, then we have that 
\beq
\label{e: |P(A,q)| exceptional}
|\ov{\mathcal{P}}(A, q)|\;=\; (2^{\omega(\kappa)-1} + 1 + \ve)\frac{h((\ve-4)\kappa^2)}{2+\ve}.
\eeq
\end{prop}

\pf\  Since $\mathcal{P}(A, q)$ is non-empty, there is a $\theta \in \mathcal{P}(A)$ with 
$q_{(A,\theta)}\sim q$. Note that $q$ is a ternary form by Theorem \ref{th:prod}.

To prove \sref{e: |P(A,q)|}, suppose first that $q$ does not represent $1$. Then \sref{e: |P(A,q)|} follows directly from Theorem 3 and formula (7) of \cite{K-CAS}. (Note that in this case 
\sref{eq:main2} reduces to the formula (3) of  \cite{K-CAS}; cf.\ Corollary 18 of \cite{K-PP}.)
  
Next, if $q$ represents $1$ but is not equivalent to $x^2 +4\kappa(y^2+\ve yz + z^2)$, for any 
$\kappa > 1$ and $\ve \in\{0,1\}$, then from the proof of Corollary 18 of \cite{K-PP} we see that
$a(q) = 2\max(1, r^*_4(q))$, and hence \sref{e: |P(A,q)|} follows from Theorem 1 of \cite{K-PP}
together with formula (7) of \cite{K-CAS}.

We are thus left with the case that $q\sim x^2+4\kappa(y^2+\ve yz+z^2)$, for some $\kappa>1$ and 
$\ve\in\{0,1\}$. Put $q_{\ve,\kappa}:=\kappa(y^2+\ve yz+z^2)$, so $\cont(q_{\ve,\kappa})=\kappa$
and $d(q_{\ve,\kappa})/\kappa^2 = \ve-4$. Since $x^2+4\kappa(y^2+\ve yz+z^2)
= f_{q_{\ve,\kappa}}$ in the notation of \cite{K-JT}, the formula given in the proof
of Proposition 53 of \cite{K-JT} shows that $I_1(q) = I_1(f_{q_{\ve,\kappa}}) = -16\kappa$ and
$d(q) = d(f_{q_{\ve,\kappa}}) = 16\kappa^2(\ve-4)$. It thus follows from Proposition \ref{p:invars}
that in this case $\De = d(q_{E,E'}) = \kappa^2(\ve-4)$ and $\De' = \ve-4$.  

 Put $u(\Delta') := |\Aut^+(q_{\ve,\kappa}/\kappa)|$ (as in \cite{K-PP}). Thus,  
$u(\Delta') = 4+2\ve$, as is well-known; cf.\ \cite{Jo}, Theorem 51a. Moreover, since $h(\De') = 1$
(cf.\ Theorem 7.30 of \cite{Cox}), we also have that $g(\De') = 1$, and so we obtain from Proposition 20 of \cite{K-PP} that 
$$
|\ov\cP(A,q)| 
= (2^{\om(\kappa)} + u(\De') - 2)\frac{2g(\De')h(\De)}{|\Aut(q_{\ve,\kappa})|h(\De')}
= (2^{\om(\kappa)} + 2+2\ve)\frac{2h(\De)}{|\Aut(q_{\ve,\kappa})|}.
$$
This shows that \sref{e: |P(A,q)| exceptional} holds because 
$|\Aut(q_{\ve,\kappa})| = 2u(\De') = 8+4\ve$. 
\spm\sps

It follows from the above proposition together with Proposition \ref{p:invars} that
$|\ov{\mathcal{P}}(A_i, q)|$ does not depend on the choice of 
$\langle A_i \rangle \in \mathcal{N}(q)$. 

\begin{cor}
\label{c:independent}
Let $A$ be a CM abelian product surface. If ${\mathcal{P}}(A, q)$ is nonempty, then  
    \beq
    \label{e: |P(a,q)| independent}
    |\ov{\mathcal{P}}(A, q)|\;=\;|\ov{\mathcal{P}}(A', q)|, \text{ for every } \langle A' \rangle \in \mathcal{N}(q).
    \eeq
\end{cor}

\pf\ By hypothesis, $A \simeq E_1\times E_2$, for some isogeneous CM elliptic curves $E_i/K$. 
Since $\mathcal{P}(A, q)\ne\emptyset$, there is a $\theta\in\mathcal{P}(A)$ such that $q_{(A,\theta)}\sim q$, and so $\langle A \rangle \in \mathcal{N}(q)$.  

Let $\langle A' \rangle \in \mathcal{N}(q)$. Then by definition there exists a 
$\theta'\in\mathcal{P}(A')$ such that $q_{(A',\theta')} \sim q \sim q_{(A,\theta)}$, and hence 
$A'$ is again a CM abelian product surface by Theorem \ref{th:prod}. Thus, $A' \simeq E'_1\times E'_2$,
for some isogeneous CM elliptic curves $E'_i/K$. Since $q_{(A',\theta')} \sim q_{(A,\theta)}$, we have that $d(q_{(A',\theta')}) = d(q_{(A,\theta)})$, $\cont(q_{(A',\theta')}) = \cont(q_{(A,\theta)})$
and that $I_1(q_{(A',\theta')}) = I_1(q_{(A,\theta)})$. It thus follows from Proposition \ref{p:invars}
that
\beq
\label{eq:|P(A,q)|}
d(q_{E'_1,E'_2})\;=\;d(q_{E_1,E_2})\quad \mbox{ and }\quad
\cont(q_{E'_1,E'_2})\;=\;\cont(q_{E_1,E_2}).
\eeq
Now if $q_{(A',\theta')}$ is not equivalent to $x^2+4\kappa(y^2+\ve yz+z^2)$, for any $\kappa > 1$ and $\ve\in\{0,1\}$, then the same is true
for $q_{(A,\theta)}\sim q_{(A',\theta')}$, and then formula \sref{e: |P(A,q)|}, together 
with \sref{eq:|P(A,q)|}, shows 
that $|\ov{\mathcal{P}}(A, q)|=|\ov{\mathcal{P}}(A', q)|$. On the other hand, if 
$q_{(A',\theta')} \sim x^2+4\kappa(y^2+\ve yz+z^2)$, for some $\kappa > 1$ 
and $\ve\in\{0,1\}$, then the assertion follows from \sref{e: |P(A,q)| exceptional}. 
\spm\sps

We are now ready to prove Theorem \ref{th:main*}.

\spm\sps

\pff\ Theorem \ref{th:main*}. By hypothesis, $A$ satisfies condition (i) of Theorem \ref{th:prod},
so by that theorem there exist two isogeneous CM elliptic curves $E_i/K$ such that 
$A\simeq E_1\times E_2$. In addition, $q := q_{(A,\theta)}$ is a ternary form. 
Thus $N_{(A,\theta)} <\infty$ by Proposition \ref{p:NAth}.  

Let $\Delta = d(q_{E_1,E_2})$ and $\kappa = \cont(q_{E_1,E_2})$. Then from \sref{eq:invars1}, 
 \sref{eq:invars2} and the definitions we see that $\Delta_q = \Delta$ and that $\kappa_q = \kappa$. 
Put $\Delta'=\Delta/\kappa^2$.

Now if $q$ is not equivalent to  $f_{q_{\ve,\kappa}}:= x^2 +4\kappa(y^2+\ve yz + z^2)$, for any $\kappa > 1$  and $\ve \in\{0,1\}$, then  \sref{eq:main*1} follows from \sref{e: |H(q)|=summation}, \sref{e: |P(a,q)| independent}, \sref{e: |A(q)|} and \sref{e: |P(A,q)|} because
\begin{align*}
   N_{(A,\theta)} &\;=\;\sum _{\langle A_i \rangle\in \mathcal{N}(q)} |\ov{\mathcal{P}}(A_i, q)|
\;=\;|\mathcal{N}(q)||\ov{\mathcal{P}}(A, q)|\\ 
&\stackrel{\text{(\ref{e: |A(q)|})}}{\;=\;}
\frac{h(\Delta')}{g(\Delta')}|\ov{\mathcal{P}}(A, q)|\stackrel{\text{(\ref{e: |P(A,q)|})}}{\;=\;}
\frac{h(\Delta')}{g(\Delta')}\frac{2^{\omega(\kappa)+1}g(\Delta')h(\Delta) a(q)}{|\Aut(q)|h(\Delta')}\\ &\;=\; 2^{\omega(\kappa)+1}h(\Delta)\frac{a(q)}{|\Aut(q)|}\;=\; 
2^{\omega(\kappa_q)+1}h(\Delta_q)\frac{a(q)}{|\Aut(q)|}.
\end{align*} %which verifies \sref{eq:main*1}.

On the other hand, if  $q \sim f_{q_{\ve,\kappa}}$, for some $\kappa>1$ and some $\ve\in\{0,1\}$, 
then $h(\De') = g(\De') = 1$, as we saw in the proof of Proposition \ref{p: |P(A,q)|}. This implies 
by \sref{e: |A(q)|} that $|{\cal N}(q)| = 1$, so $N_{(A,\theta)} = |\ov\cP(A, q)|$ 
by \sref{e: |H(q)|=summation}, and hence
\sref{eq:main*2} follows from \sref{e: |P(A,q)| exceptional}.   
\spm\sps

We will now deduce Theorem \ref{th:main} from Theorem \ref{th:main*}. For this, recall (from \cite{Mi}, for example) that if $C/K$ is a curve of genus $2$, then its Jacobian $J_C$ is an abelian surface
and the image $j(C)$ of $C$ via the canonical embedding $j: C \hookrightarrow J_C$ is an ample divisor on $J_C$. Furthermore, the image $\theta_C \in \NS(J_C)$ of $j(C)$ in the N\'eron-Severi group is a principal polarization. We then write $q_C:= q_{(J_C,\theta_C)}$ for its associated refined Humbert invariant.  

It turns out that whether or not a given principally polarized abelian surface $(A,\theta)$ is a Jacobian can be determined from its refined Humbert invariant $q_{(A,\theta)}$. Indeed, by
by Proposition 6 of \cite{K-MJ} we have that 
\begin{equation} 
\label{e: theta irreducible}
(A,\theta) \simeq (J_C,\theta_C), \mbox{ for some curve }C/K\; \Leftrightarrow\; 
q_{(A,\theta)}(D)\neq 1, \forall D\in \NS(A,\theta).
 \end{equation}   

This result, together with Torelli's Theorem, implies that the number $N_C$ defined in the introduction equals the number $N_{(J_C,\theta_C)}$, if $C$ is of CM product type. More precisely:   

\begin{lemma}
\label{lem:NC}
If $C/K$ is a curve of genus $2$, then the map $C' \mapsto (J_{C'},\theta_{C'})$ induces a bijection
${\cal N}_C \stackrel\sim\ra {\cal N}_{(J_C,\theta_C)}$, where $\mathcal{N}_C$ denotes the set of isomorphism classes $\langle C' \rangle$ of curves $C'/K$ of genus $2$ such that $q_{C'}\sim q_C$. 
 Thus, if $C/K$ is of CM product type, then 
$N_C := |{\cal N}_C| = N_{(J_C,\theta_C)}$.
\end{lemma}

\pf\ It is clear from the definitions that the given rule induces a map 
$\vp_C: {\cal N}_C \ra {\cal N}_{(J_C,\theta_C)}$. This map is injective by Torelli's Theorem; cf.\
Theorem 12.1 in \cite{Mi}. Moreover, $\vp_C$ is surjective because if 
 $\langle A,\theta \rangle \in \mathcal{N}_{(J_C,\theta_C)}$, then 
$q_{(A,\theta)} \sim q_{(J_C,\theta_C)}$. Since $q_{(J_C,\theta_C)}$ does not represent 1 
by \sref{e: theta irreducible}, it follows that also $q_{(A,\theta)}$ does not represent 1, 
and so by  \sref{e: theta irreducible} again we have that $(A,\theta) \simeq (J_{C'},\theta_{C'})$,
for some curve $C'/K$ of genus $2$. Since $q_{C'} = q_{(J_{C'}, \theta_{C'})} \sim 
q_{(A,\theta)}\sim q_{(J_C,\theta_C)} = q_C$, we see that $\lb C'\rb \in {\cal N}_C$, so 
$\lb A,\theta\rb = \vp_C(\lb C'\rb)$ lies in the image of $\vp_C$. Thus, $\vp_C$ is surjective and hence bijective. 

If $C$ is of CM product type, then $J_C$ is a CM abelian product surface by Theorem \ref{th:prod}, and so $N_{(J_C,\theta_C)} < \infty$ by Proposition \ref{p:NAth}. Thus, 
the first assertion shows that $N_C = |{\cal N}_{(J_C,\theta_C)}| =
N_{(J_C,\theta_C)}$.  

\spm\sps

\pff\ Theorem \ref{th:main}. Since $C$ is of CM product type, we have by Lemma \ref{lem:NC} 
that $N_C = N_{(J_C,\theta_C)}$. As was mentioned in its proof, $J_C$ is a CM abelian product surface. 
By \sref{e: theta irreducible} we know that $q_C = q_{(J_C,\theta_C)}$ cannot represent $1$, so $q_C$ cannot be equivalent to one of the exceptional forms of Theorem \ref{th:main*}, and hence it follows from that theorem that  
$$
N_C \;=\; N_{(J_C,\theta_C)} \;=\;  
2^{\omega(\kappa_{q_C}) + 1}h(\Delta_{q_C})\frac{a(q_C)}{|\Aut(q_C)|}.
$$
Since $q_C$ represents a square by Theorem \ref{th:prod}, we have by Theorem 25 of 
\cite{K-CAS} that 
\beq
\label{eq:aq-Aut}
2a(q_C) \;=\; |\Aut(C)|,
\eeq
and so \sref{eq:main} follows. 
\spm\sps

\pff\ Corollary \ref{c:main}. To prove the first assertion, fix an integer $n\ge1$ and let 
$$
Q_n \;=\; \{q_C : C/K \mbox{ is a curve of genus $2$ of product CM type with $N_C \le n$}\}/\!\sim  
$$
denote the set of equivalence classes of those refined Humbert invariants which arise from curves 
$C/K$ of genus $2$ of product CM type with $N_C \le n$. It clearly suffices to show that 
$|Q_n|<\infty$ because by definition the number of isomorphism classes of such curves with a given
$q\in Q_n$ is equal to $N_C \le n$.

Now if $q\in Q_n$, then by Theorem \ref{th:prod} we know that $q$ is a positive ternary form 
and so $|\Aut(q)|\le 48$ because $|\Aut^+(q)|\leq 24$ by Theorem 105 of \cite{Di} and because $|\Aut(q)| = 2|\Aut^+(q)|$ since $-1\in \Aut(q)\setminus\Aut^+(q)$. 
We thus obtain from \sref{eq:main} that if $q$ is in $Q_n$, then 
$$
\frac{h(\Delta_q)}{48}\; \le \; 2^{\omega(\kappa_q)}h(\Delta_q)\frac{|\Aut(C)|}{|\Aut(q)|} 
\;=\; N_C \; \le \; n.
$$
Now if $q\in Q_n$, then $\Delta_q = d(q)/16$ is a negative quadratic discriminant by 
Proposition \ref{p:invars}. But by the celebrated result of Heilbronn \cite{H} we know that there are only finitely many negative discriminants $\Delta$ such that $h(\Delta)\leq 48n$, so there exists
a bound $B_n$ such that $h(\Delta)\leq 48n \Ra |\De| \le B_n$. We therefore see that $Q_n$ is a subset of the set
$$
\tilde Q_n = \{q: q \mbox{ is a positive ternary form with $|d(q)| \le 16B_n$}\}/\!\sim.
$$     
Now by Theorem 11 of \cite{Wa},  there are up to equivalence only finitely many forms of a given rank and discriminant, so $\tilde Q_n$ is a finite set, and hence so is its subset $Q_n$. This proves
the first assertion.   

In order to prove the second assertion, it suffices in view of the first assertion to show that there exist infinitely many non-isomorphic genus 2 curves $C/K$ of CM product type. 

To verify this, we first show that the set $\cA_{CM}(K) = \{\lb A\rb\}$ consisting of the set of   
isomorphism classes of CM abelian product surfaces $A/K$ is infinite. 
If $\Char(K) = 0$, then this follows immediately from Theorem 71 of \cite{K-PCM} (and its proof) because there are clearly infinitely many negative discriminants $\De \equiv 0,1\Mod 4$, and hence there are infinitely many equivalence classes of positive binary quadratic forms. If $\Char(K) = p > 0$, 
then this follows from Remark 73 of \cite{K-PCM} because there exist infinitely many negative numbers $\De$ with $\De \equiv 1\Mod{4p}$, so there exist infinitely many positive binary quadratic forms whose discriminant $\De$ satisfies $(\frac{\De}p) = 1$.  
    
Next, let $\cA_{CM}^*(K)$ denote the subset of $\cA_{CM}(K)$ consisting of the isomorphism classes of those CM abelian product surfaces $A/K$ which contain a smooth genus $2$ curve $C/K$. 
By Theorem 2 of \cite{K-JT}
we know that $|\cA_{CM}(K) \setminus \cA^*_{CM}(K)| \le 15$, so $\cA_{CM}^*(K)$ is also an infinite set.  

Now if $C/K$ is a smooth genus $2$ curve lying on an abelian surface $A/K$, then $A \simeq J_C$ 
(by the universal property of the Jacobian), and so if $C \subset A$, where $\lb A\rb\in\cA^*_{CM}(K)$, then $C$ is a curve of CM product type. Furthermore, since $C\simeq C'$ implies that 
$J_C \simeq J_{C'}$, it therefore follows that the infinitely many non-isomorphic 
abelian surfaces in $\cA^*_{CM}(K)$ give rise to infinitely many non-isomorphic curves $C/K$ 
of CM product type, as claimed. This proves the second assertion. 

%**********************

\section{Examples}
\label{s:examples}

In this section, we want to illustrate how Theorem \ref{th:main} can be applied to some explicit cases.
%We begin with the following result which will be used below.
In this context, the following result discusses some cases satisfying $N_C =1$.  

\begin{prop}
\label{p:example}
Let $A/K$ be a CM abelian product surface, and let $\theta\in\cP(A)$ be a principal polarization on $A$.
Suppose that $\kappa_{q_{(A,\theta)}} = 1$ and that $h(\De) = 1$, where $\De = \De_{q_{(A,\theta)}}$. 
Then $A \simeq E\times E$, where $E/K$ is up to isomorphism the unique CM elliptic curve such that 
$\End(E)$ has discriminant $\De$.  

Conversely, if $E/K$ is a CM elliptic curve such that $h(\De) = 1$, where $\De$ is the discriminant of 
$\End(E)$, then $A = E\times E$ is a CM abelian product surface, and $\De_{q_{(A,\theta)}} = \De$,
$\kappa_{q_{(A,\theta)}} = 1$, for every principal polarization $\theta \in\cP(A)$ on $A$. In addition, we have that $N_{(A,\theta)} = 1$ and $|\Aut(q_{(A,\theta)})| = 2a(q_{(A,\theta)})$. 

In particular, if $C$ is a curve of CM product type with $\kappa_C = 1$ and $h(\De_C)=1$, then 
$N_C =1$, $|\Aut(C)| = |\Aut(q_C)|$ and $J_C \simeq E\times E$, where $E/K$ is up to isomorphism the unique elliptic curve such that $\End(E)$ has discriminant $\De_C$.  
\end{prop}

\pf\ Since $A/K$ is a CM abelian product surface, we have that $A\simeq E_1\times E_2$, where $E_i/K$ 
are two isogeneous CM elliptic curves $E_i/K$. By Proposition \ref{p:invars} we have that
$d(q_{E_1,E_2}) = d(q_{(A,\theta)})/16 = \De_{q_{(A,\theta)}} = \De$ and $\cont(q_{E_1,E_2}) 
= \kappa_{q_{(A,\theta)}} = 1$. Thus, $q_{E_1,E_2}$ is a primitive positive binary quadratic form of
discriminant $\De$. Since $h(\De)=1$ by hypothesis, $q_{E_1,E_2}$ is equivalent to the principal
form $1_\De$ of discriminant $\De$, and hence $q_{E_1,E_2}$ represents $1$. Thus, there exists
an $h\in\Hom(E_1,E_2)$ with $\deg(h)=1$, so $E:=E_1\simeq E_2$, and hence $q_{E_1,E_2} \sim q_{E,E}$. This means that $\End(E)$ has discriminant $\De$, and that $A\simeq E\times E$.

Since $h(\De) = 1$, then, as is well-known, the isomorphism class of the elliptic curve $E/K$ is uniquely determined by the fact that the endomorphism ring of $E/K$ has discriminant $\De$; cf.\
\cite{K-PCM}, equation (55). 

Conversely, suppose that $A = E\times E$ where $E/K$ is a CM elliptic curve, and let $\De$ be the discriminant of $\End(E)$. Then $\De = d(q_{E,E})$ (by definition), and clearly $\cont(q_{E,E}) =1$
because $1_E \in\End(E)$. Thus, by Proposition \ref{p:invars} (and the definitions) we have that
$\De_{q_{(A,\theta)}} = \De$ and $\kappa_{q_{(A,\theta)}} = 1$, for every $\theta\in \cP(A)$. 
   
It remains to show that $N_{(A,\theta)} = 1$ and that $2a(q) = |\Aut(q)|$, where $q= q_{(A,\theta)}$. 
Now since $\kappa_q =1$, we see that $q$ is not one of the exceptional forms of 
Theorem \ref{th:main*}, and so it follows from \sref{eq:main*1} that 
$N_{(A,\theta)} = \frac{2a(q)}{|\Aut(q)|}$. Since $a(q) = a(\theta)$  by Corollary 17 of \cite{K-PP} (where $a(\theta)$ is defined as in \cite{K-PP}) and 
since $a(\theta)\mid |\Aut^+(q)| = \frac12|\Aut(q)|$ by Proposition 16 of \cite{K-CAS}, we see that
$N_{(A,\theta)} \le 1$. But since $N_{(A,\theta)} \ge 1$ (because 
$\lb A,\theta\rb \in \cN_{(A,\theta)}$), it follows that $N_{(A,\theta)} = 1$ and 
$2a(q) = |\Aut(q)|$, as claimed.

The last assertion follows directly from the first part by using equation \sref{eq:aq-Aut}.
%the fact that $2a(q_C) = |\Aut(C)|$; cf.\ Theorem 25 of \cite{K-CAS}.     
  
\begin{rmk} \em
\label{r:h=1}
(a) If $\Char(K)=0$ and if $E/K$ is a CM elliptic curve such that $h(\De_E) = 1$, where $\De_E$ is the discriminant $\End(E)$, then $E/K$ comes from an elliptic curve defined over $\Q$. An equation of such a curve, together with its $j$-invariant, is given in the tables on p.\ 483 of Silverman \cite{Si}.

Similarly, if $\Char(K) = p>0$, then the CM elliptic curves with $h(\De_E) = 1$ are obtained by reduction mod $p$ from the corresponding CM curve
over $\Q$, provided that $(\frac{\De_E}p) = 1$. On the other hand, if $(\frac{\De_E}p) \ne 1$, then there is no such CM elliptic curve.
\spm

(b) In the situation of Proposition \ref{p:example} there may not be a curve $C/K$ satisfying the given conditions. Indeed, if $\De = -3,-4$ or $-7$, then there is no such curve; cf.\ 
Hayashida/Nishi \cite{HN-cm}, Theorem 1, or \cite{K-JT}, Theorem 2, or the table in \cite{K-CAS}.    
\end{rmk}

We now want to show how Proposition \ref{p:example} can be used to deduce the uniqueness of a curve 
$C/K$ when we start with a fixed ternary quadratic form.

\begin{exam}
\label{ex:1}
Consider the following ternary quadratic form 
$$
q(x,y,z) = 4x^2+4y^2+4z^2+4yz+4xz+4xy.
$$
If $\Char(K)=0$, then up to isomorphism there is a unique genus $2$ curve $C/K$ such that $q_C \sim q$.
Moreover, $|\Aut(C)| = 48$ and the Jacobian $J_C$ of $C$ is isomorphic to $E \times E$, where $E$ is given by the equation $y^2 = x^3 + 4x^2 + 2x$ and has $j$-invariant $j(E) = 2^6\cdot3^5$.   
\end{exam}

\pf\  We first note that $q/2$ is an \emph{improperly primitive} form (in the sense of 
Dickson \cite{Di}) because $q/4$ is a primitive form and the coefficient of the term $yz$ (or of $xy$, 
or of $xz$) of $q/4$ is odd. Thus, $q$ satisfies condition (1.1) in Theorem 1.2 of \cite{Ki}, and also 
condition (1.2) holds because $q(1,0,0) = 4$ is a square. Since $\Char(K) =0$, it thus follows from 
that theorem that there exists a principally polarized abelian surface $(A,\theta)$ over $K$ such 
that $q_{(A,\theta)} \sim q$. 

It is clear that $q$ cannot represent $1$, so it follows from  \sref{e: theta irreducible} that 
$(A,\theta) \simeq (J_C,\theta_C)$, for some genus $2$ curve $C/K$. Note that since $A/K$ is a CM abelian product surface by Theorem \ref{th:prod}, it follows that $C$ is a curve of CM product type, 
so $J_C \simeq E_1\times E_2$, for some isogeneous CM elliptic curves $E_i/K$.

In order to apply Proposition \ref{p:example}, we need to show that $C$ satisfies the required
conditions. For this, we first compute the discriminant $d(q)$. By definition (cf.\ \cite{Wa}, p.\ 2)
we have that $d(q) = -\frac12\det(M(q))$, where $M(q)$ is the coefficient matrix of $q$. In our case 
we have that 
$$
M(q) \;=\; \begin{pmatrix}
8 & 4 & 4\\
4 & 8 & 4\\
4 & 4 & 8
\end{pmatrix},
$$
so we see that $d(q) = -8\cdot16$. Thus, by Proposition \ref{p:invars} we obtain that $d(q_{E_1,E_2}) = -8$.
Since $-8$ is a fundamental discriminant, it follows that $q_{E_1,E_2}$ is primitive, so
$\kappa_q = \cont(q_{E_1,E_2}) = 1$. Thus, $C$ satisfies the hypotheses of Proposition 
\ref{p:example} because $h(-8) = 1$ (cf.\ Theorem 7.30 of \cite{Cox})), and so $J_C \simeq E \times E$,
where $E$ is a CM elliptic curve such that $\End(E)$ has discriminant $-8$. 
%(so $\End(E) \simeq \Z[\sqrt{-2}]$). 
Since this determines $E/K$ uniquely up to isomorphism (cf.\ Proposition \ref{p:example}), we see
from the first table on p.\ 483 of \cite{Si} that $j(E) = 2^6\cdot5^3$ and from the second table there that
$E/K$ is given by the equation $y^2 = x^3+4x^2+2x$.   

To determine $\Aut(C)$, we can use either of two methods. The first consists of computing 
$r^*_4(q)$. For this, we observe that clearly 
$$\{\pm(1,0,0),\pm(0,1,0), \pm(0,0,1), \pm(1,-1,0), \pm(1,0,-1), 
\pm(0,1,-1)\} \subset R_4(q),
$$ so $r^*_4(q) \ge 12$. We thus see from \sref{eq:main2} 
that $a(q) = 3r^*_4 -12 \ge 24$, so by \sref{eq:aq-Aut} we obtain that $|\Aut(C)| = 2a(q) \ge 48$. 
But since $|\Aut(C)| \le 48$, for all
such curves $C/K$ (cf.\ Theorem 25 of \cite{K-CAS}), it follows that $|\Aut(C)| = 48$ (and that 
$r^*_4(q) = 12$).     
  
Alternately, we can determine $|\Aut(C)|$ by calculating $|\Aut(q)| = |\Aut(q/2)|$, 
and for this we can use the method and tables of Dickson \cite{Di}, \S82-83. Indeed, we see that
$q/2$ is the first entry of Table II on p.\ 185 of \cite{Di}, so we have by that table that
$|\Aut^+(q/2)| = 24$. (The method for computing this table is given in Theorem 105 of \cite{Di}.)
We thus have that $|\Aut(q)| = 2|\Aut^+(q/2)| = 48$, and hence also $|\Aut(C)| = 48$ by Proposition
\ref{p:example}.
  
\begin{rmk} \em
\label{r:Burnside}   
One can say more about the curve $C/K$ which was constructed in Example \ref{ex:1}. Indeed, 
since here $\Char(K) =0$,
it is known that up to isomorphism there is only one curve $C/K$ with $|\Aut(C)| = 48$, and so
it follows that $C$ is given by the equation $y^2 = x(x^4-1)$; cf.\ \cite{Ig}, \S8. 
 Furthermore, we have that $\Aut(C) \simeq \GL_2(3)$ by Theorem 2 of \cite{SV}. In \cite{AP}, this 
curve is called the \emph{Burnside curve}.
\end{rmk}

\iffalse
As was discussed in Example \ref{ex:1}, we deduce from Proposition \ref{p:example}  the uniqueness of a curve $C/K$ for given  ternary quadratic forms in the following example. This example also illustrates that there may be two non-isomorphic curves of genus 2 on the same abelian surface.
\fi

The next Example \ref{ex:Delta=28} uses again Proposition \ref{p:example} to deduce the uniqueness of a curve $C/K$ for a given ternary quadratic form. This example also illustrates that there may be two non-isomorphic curves of genus 2 on the same abelian surface. 

\begin{exam}
\label{ex:Delta=28}
    Consider the following two ternary quadratic forms
$$
q_1(x,y,z) \;=\; 4x^2+5y^2+8z^2-4yz-4xy \text{ and } q_2(x,y,z) \;=\; 4x^2+4y^2+8z^2-4xz.
$$    
  If $\Char(K) = 0$, then up to isomorphism there is a unique genus $2$ curve $C_{i}/K$ such that $q_{C_i}\sim q_i$, for $i=1,2$. Moreover, $|\Aut(C_1)| = 4$ and $|\Aut(C_2)| = 8$, and the Jacobian $J_{C_i}$ of $C_i$ is isomorphic to $E\times E$, where $E$ is given by the equation $y^2 = x^3 - 595x + 5586$ and has $j$-invariant $j(E) = 3^3\cdot5^3\cdot17^3$. Furthermore, $C_1\not\simeq  C_2$.
\end{exam}

\pf\ For the first form $q_1$, note first that $q_1(x,y,z) = (2x - y)^2 + 4(y^2 - yz + 2z^2)$ is the sum of two positive forms, and so $q_1$ is a positive form. In fact, we easily see that $q_1\geq 4$, if $(x,y,z) \neq (0,0,0)$, so $q_1$ cannot represent 1.

As in Example \ref{ex:1}, we can calculate the discriminants $d(q_i)$, for $i=1, 2$ by computing the determinant of $M(q_i)$. This gives $d(q_i) = -16\cdot 28$, for $i=1 ,2$. Since $q_1(1,1,1) = 3^2$, the form $q_1$ represents a square which is relatively prime to its discriminant. 
Furthermore, since $q_1 \equiv y^2 \Mod{4}$, we see that $q_1 \equiv 0, 1 \Mod{4}$. 
If $\Char(K) =0$, then it follows from Theorem 1 of \cite{K-rH} that $q_1 \sim q_{(A,\theta)}$, 
for some principally polarized abelian surface $(A,\theta)$ over $K$. Moreover, as was mentioned above, $q_1$ cannot represent 1, so it follows from \sref{e: theta irreducible} that 
$(A,\theta) \simeq (J_{C_1}, \theta_{C_1})$, for some genus 2 curve $C_1/K$. 
As in Example \ref{ex:1}, we conclude that $C_1$ is of CM product type.

For the second form $q_2$, note that $q_2/2$ is an improperly primitive form because $q_2/4$ is a primitive form and the coefficient of the term $xz$ of $q_2/4$ is odd. Moreover, $q_2$ is a positive form because $q_2(x,y,z) = 4(x^2 -xz + 2z^2) + 4y^2$ is the sum of two positive forms. 
Thus, $q_2$ satisfies condition (1.1) in Theorem 1.2 of \cite{Ki}, and also condition (1.2) holds because $q_2(1,0,0) = 2^2$. If $\Char(K) = 0$, then it follows from that theorem that there exists a principally polarized abelian surface 
$(A',\theta')$ over $K$ such that  $q_2 \sim q_{(A',\theta')}$. Furthermore, it is clear that $q_2$ cannot represent 1, so it again follows from  \sref{e: theta irreducible} that 
$(A',\theta') \simeq (J_{C_2}, \theta_{C_2})$, for some genus 2 curve $C_2/K$. Similarly, we conclude that $C_2$ is of CM product type.

In order to apply Proposition \ref{p:example}, we need to show that the $C_i$ satisfy the required conditions.

To begin, we have that $\Delta_{C_i} = d(q_i)/16 = -28$, and so $h(\Delta_{C_i}) = 1$ by Theorem 7.30 of \cite{Cox}. Next, to compute $\kappa_{C_i}$, recall from \S\ref{s:forms} that the genus invariant $|I_1(q)|$ of a primitive form $q$ is the content of the adjoint form adj$(q)$ (see \sref{eq: adj_form}). Hence, we see that $|I_1(q_1)| = \gcd(16\cdot9,16\cdot8,16\cdot4,16\cdot4,16\cdot2,16\cdot8) = 16$ 
and that $|I_1(q_2/4)| = \gcd(8,7,4,0,4,0) = 1$, and thus $\kappa_{C_i} = 1$, for $i= 1, 2$ (by definition). This shows that the $C_i$ satisfy the hypotheses of Proposition \ref{p:example}, and so 
$J_{C_i} \simeq E_i\times E_i$, where $E_i$ is a CM elliptic curve such that $\End(E_i)$ has discriminant $-28$, for $i=1, 2$. Since this determines $E_i/K$ uniquely up to isomorphism (cf.\ Proposition \ref{p:example}), we obtain that $E_1\simeq E_2 := E$, and hence $J_{C_i}\simeq E\times E$, for $i=1, 2$. Furthermore, we see from the first table on p.\ 483 of \cite{Si} that $j(E) = 3^3\cdot5^3\cdot17^3$ and from the second table there that $E/K$ is given by the equation 
$y^2 = x^3 - 595x + 5586$, as asserted.

To determine $|\Aut(C_i)|$, we can use two different methods as was mentioned in Example \ref{ex:1}. 
We use the first method to find $|\Aut(C_1)|$, and the second one to calculate $|\Aut(C_2)|$. 

The first method consists of computing $a(q_1)$, and this can be achieved by calculating $r_4^*(q_1)$. But, this is already done in Corollary 30 of \cite{K-CAS}, for primitive ternary forms $q_C$. Indeed, since $q_1$ satisfies the inequalities of Theorem 103 of \cite{Di}, $q_1$ is an 
\emph{Eisenstein-reduced} ternary form. Since $q_1$ is primitive, we are in the situation of 
Corollary 30 of \cite{K-CAS}, and it follows from there that $a(q_1) = 2$. Thus, 
$|\Aut(C_1)| = 2a(q_1) = 4$ by equation \sref{eq:aq-Aut}.

To apply the second method to the form $q_2$,  we determine $|\Aut(q_2)| = |\Aut(q_2/2)|$. %To this end, we can use the method and tables of Dickson\cite{Di}, \S.82-83. Indeed,
Note that $q_2/2$ is the sixth entry of Table II on p.\ 185 of \cite{Di}, and so we have by that table that $|\Aut^+(q_2/2)| = 4$, and thus $|\Aut(q_2)| = 2|\Aut^+(q_2/2)| = 8$. Therefore, it follows that 
$|\Aut(C_2)| = |\Aut(q_2)| = 8$ by Proposition \ref{p:example}. 

Finally, we observe that $C_1 \not\simeq C_2$ because $q_{C_1} \not\sim q_{C_2}$. Indeed, $q_{C_1} \sim q_1$ is primitive, whereas $q_{C_2} \sim q_2$ is not, so $q_{C_1} \not\sim q_{C_2}$.
\spm

The following example shows that it is possible to have $N_C = 1$ in a situation where
Proposition \ref{p:example} does not apply. %while $|\Aut(C)|\neq |\Aut(q_C)|$ as follows.

\begin{exam} 
\label{ex:2}
Let $q$ be the primitive ternary quadratic form defined by 
$$
q(x,y,z) \;=\; 4x^2+4y^2+9z^2-4xy.
$$
If $\Char(K) = 0$ or if $\Char(K) \equiv 1\Mod 3$, then there exists (up to isomorphism) a unique curve 
$C/K$ such that $q_C \sim q$. Furthermore, $|\Aut(C)| = 12$  and $J_C \simeq E_1\times E_2$,
where $E_1$ is the elliptic curve given by $y^2 + y = x^3$ and $E_2$ is the elliptic curve given by
$y^2 + y = x^3 -30x + 63$. Moreover, $j(E_1) =0$ and $j(E_2) = -2^{15}\cdot3\cdot 5^3$.   
\end{exam}

\pf\ Note first that $q(x,y,z) = 4(x^2 -xy + y^2) + 9z^2$ is the sum of two positive forms, so $q$ is also positive. In fact, we see that $q(x,y,z) \ge 4$, if $(x,y,z)\ne (0,0,0)$, so $q$ cannot represent 
$1$.

 As in Example \ref{ex:1}, we can calculate its discriminant by computing the determinant of 
$M(q)$. This gives $d(q) = -16\cdot 27$.  Since $q(2,2,1) = 5^2$, the form $q$ represents a square which is relatively prime to its discriminant. Furthermore, since $q \equiv z^2 \Mod 4$, we see that
$q\equiv 0, 1\Mod 4$. If $\Char(K) =0$, then it follows from Theorem 1 of \cite{K-rH} that 
$q \sim q_{(A,\theta)}$, for some principally polarized abelian surface $(A,\theta)$ over $K$.
Similarly, if $p := \Char(K) >0$, then $(\frac{d(q)/16}{p}) = (\frac{-3}p) = 1$ by our hypothesis, 
and so it follows from Theorem 28 of \cite{K-rH} that $q \sim q_{(A,\theta)}$, for some principally polarized abelian surface $(A,\theta)$ over $K$.

As was mentioned above, $q$ cannot represent $1$, so it follows from \sref{e: theta irreducible} 
that $(A,\theta) \simeq (J_C,\theta_C)$, for some genus $2$ curve $C/K$. As in Example \ref{ex:1},
we conclude that $C$ is of CM product type.

To show that $C/K$ is uniquely determined by $q$, it suffices to show that
$N_C =1$, and for this we will use the formula \sref{eq:main}. Thus, we need to compute the quantities on the right hand side of \sref{eq:main}.  

To begin, we have that $\De_C = d(q)/16 = -27$, and so $h(\De_C) = 1$ (by Theorem 7.30 of \cite{Cox}). 
Next, we see that $|I_1(q)| = \gcd(16\cdot9,16\cdot9,48,0,0,16\cdot9) = 48$ as was discussed in Example \ref{ex:Delta=28}, and thus $\kappa_C = 3$.
  
It remains to compute $|\Aut(C)|$ and $|\Aut(q)|$.  For this, observe that $q$ is an Eisenstein-reduced primitive ternary form, and we are thus in the situation of Corollary 30 of \cite{K-CAS}, and so it
follows from there that $a(q) = 6$. Thus, $|\Aut(C)| = 2a(q) = 12$ by equation \sref{eq:aq-Aut}.

Next, to compute $|\Aut(q)|$, we will use Theorem 105 of \cite{Di}. Since $q$ is Eisenstein-reduced,
 we see that the table on p.\ 180 of \cite{Di} is applicable. We observe that the cases of lines 2, 5 and 9 hold and that the others do not. 
These cases give us five automorphs, and hence we have six automorphs including the identity. Since the coefficients of the terms $yz$ and $xz$ are $0$, we conclude that 
$|\Aut^+(q)| = 6\cdot2 = 12$ by the footnote on p.\ 180 of \cite{Di}, and so $|\Aut(q)| = 24$.

It therefore follows from \sref{eq:main} that 
$N_C = 2^{\om(3)}h(-27)\frac{|\Aut(C)|}{|\Aut(q_C)|} = \frac{2\cdot 12}{24} = 1$, so $C/K$ is uniquely determined by $q$ (up to isomorphism).

It remains to determine the structure of $J_C$. Since $C$ is of CM product type, we have that 
$J_C\simeq E_1\times E_2$, for some isogeneous CM elliptic curves $E_1$ and $E_2$. 
By \sref{eq:invars1} and \sref{eq:invars2} and the above computations, 
we see that $q_{E_1,E_2}$ is a binary quadratic form of content $3$ and discriminant $-27$. 
%by using the values of $d(q)$ and $I_1(q)$. 
Moreover, by \sref{eq:Nq=NA1} we have that 
$\End^0(E_1)\simeq \End^0(E_2)= \Q(\sqrt{-27}) =:F$, and so $\Delta_F = -3$. 

As in the proof of Proposition \ref{p: |A(q)|}, let $f_i = f_{E_i}$ denote the conductor of the order 
$\End(E_i)$ of $F$. Then by the first formula of (79) in \cite{K-PCM} we see that 
$\lcm(f_1,f_2)^2 = \frac{d(q_{E_1,E_2})}{\De_F} = \frac{-27}{-3} = 9$, and so it follows from the second formula of (79) in \cite{K-PCM} that 
$\gcd(f_1,f_2) = \frac{\lcm(f_1,f_2)}{\cont(q_{E_1, E_2})} = \frac33 = 1$. Thus, by interchanging 
$E_1$ and $E_2$, if necessary, we may assume without loss of generality that $
f_{1} = 1$ and $f_{2} = 3$. This means that $\De_{E_1} = -3$ and $\De_{E_2} =-27$. Since 
$h(-3) = h(-27) = 1$, the curves $E_i/K$ are uniquely determined (up to isomorphism) by their endomorphism discriminants, and so we see from the first table of \cite{Si} on p.\ 483 that $j(E_1) =0$ and that $j(E_2) = -2^{15}\cdot3\cdot 5^3$, and from the second table that $E_1$ and $E_2$ are given by the asserted equations.

%***************************biblio

%\bibliographystyle{cdraifplain}

%\bibliographystyle{elsarticle-num}
%\bibliographystyle{cas-model2-names}
\bibliographystyle{abbrvnat}
\bibliography{card}

\end{document}